\documentclass{article}

\usepackage{amsmath}
\usepackage{amssymb}
\usepackage[all]{xy}

\numberwithin{equation}{section}

\begin{document}

\centerline{\bf TOWARDS MOTIVIC QUANTUM COHOMOLOGY OF $\overline{M}_{0,S}$}

\bigskip

\centerline{\bf Yuri I.~Manin, Maxim Smirnov}

\medskip

\centerline{\it Max--Planck--Institut f\"ur Mathematik, Bonn, Germany}
\bigskip
\hfill{\it Dedicated to Professor V.~V.~Shokurov, on his 60th birthday}
\bigskip

{ABSTRACT.}  We explicitly calculate some Gromov--Witten correspondences determined by maps of labeled curves of genus zero to the moduli spaces of labeled curves of genus zero. We consider these calculations as the first step towards studying the self--referential nature of motivic quantum cohomology.

\bigskip

\centerline{\bf 0. Introduction: Motives and Quantum Cohomology}

\medskip

{\bf 0.1. A summary.} Let $Var_k$ be the category of smooth complete algebraic varieties defined over a field $k$.

\smallskip

The category of classical motives  $Mot_k^K$, with coefficients in a $\textbf{Q}$--algebra $K$, is the target of a functor $h:\,Var_k^{op}\to Mot_k^K$ which, in the vision of Alexandre Grothendieck, ought to be {\it a universal cohomology theory,} with values in a tensor $K$--linear category.

\smallskip

Morphisms $X\to Y$ in $Mot_k^K$ are represented by {\it classes of correspondences}, algebraic  cycles on $X\times Y$ with coefficients in $K$. Depending on the equivalence relation, imposed on these cycles, one could consider Chow motives, numerical motives, etc.

\smallskip

Besides objects $h(V)$ for $V\in Var_k$,  the category $Mot_k^K$ contains their  {\it direct summands}, (``pieces'') and their {\it twists} by Tate's motive.  Formally adding these objects one turns the category of motives into a Tannakian category. One can then apply to it the philosophy of motivic fundamental group. Ideally, all inherent structures of cohomology objects can be encoded/replaced by the representations of  the respective motivic fundamental group. 

\smallskip

What is special about  ``total motives''  $h(V)$, as opposed to pieces and twists? 

\smallskip

For example,  $h(V)$'s bring with them a  natural structure of {\it commutative algebras} in $Mot_k^K$. It is not determined only by $h(V)$: isomorphic motives $h(V)$'s may well have different multiplications; but of course, this classical multiplication is motivic in the sense that it is induced by the diagonal map $V\to V\times V$ in $Var_k$ and by the class of its graph in $Mot_k^K$. 

\smallskip

The advent of quantum cohomology from physics to algebraic geometry opened our eyes to the fact that classical cohomology spaces  of algebraic varieties, say, over $\textbf{C}$, possess an incomparably richer structure: they (or rather their tensor powers) are acted upon by cohomology of moduli spaces of pointed curves $H^*(\overline{M}_{g,n})$, much as Steenrod operations act in topological situation. From the physical perspective, these operations encode  ``quantum corrections to the classical multiplication''.

\smallskip

Grothendieck's vision however turned out  prophetic : {\it this new structure is  motivic as well} in the same sense: it is induced by canonical Chow correspondences, Gromov--Witten invariants $I_{g,n,\beta}^V$ in  $A_*(\overline{M}_{g,n}\times V^{n})$ indexed by effective classes $\beta$ in $A_1(V)$. This was conjectured in [KoMa1], worked out in more detail in [BehMa], and finally proved in [Beh1], where the virtual fundamental classes in the Chow groups of spaces of stable maps were constructed  by algebraic--geometric techniques.

\smallskip

This construction allowed K.~Behrend to establish a list of universal identities between the Gromov--Witten invariants that were conjectured earlier.

\smallskip

Taken together, these identities imply that for each total motive $h(V)$,   the infinite sum of its copies indexed by the numerical classes $\beta$ of effective curves on $V$ possesses the canonical structure of an algebra over the cyclic modular operad  $ \mathcal{H}\mathcal{M}$:
$$
\mathcal{H}\mathcal{M}(n):=\coprod_{g} h(\overline{M}_{g,n})
$$
This is the motivic core of quantum cohomology.

\smallskip

However, this discovery also stressed an inherent tension between the initial Grothendieck vision and the highly non--Tannakian character of the quantum cohomology  expressed in the following observations. 

\smallskip

First, these structures of $ \mathcal{H}\mathcal{M}$--algebras are not functorial in any naive sense wrt  morphisms in $Var_k$ (except isomorphisms). Notice that the classical multiplication {\it is} functorial wrt morphisms in $Var_k$; quantum corrections destroy this.  However, as was shown in [LeLW], quantum multiplication {\it is functorial} wrt at least certain isomorphisms in  $Mot_k^K$ (flops) that do not agree with classical multiplication: quantum corrections exactly compensate classical discrepancies. This is a remarkable fact suggesting that motivic functoriality might be an important hidden phenomenon.

\smallskip

Second,  being total motives, $h(\overline{M}_{g,n})$ {\it themselves} have quantum cohomology, that is, define algebras over $\mathcal{H}\mathcal{M}$.

\smallskip

One aim of this note is to draw attention to this self--referentiality and  to start studying the quantum cohomology {\it of} $\mathcal{H}\mathcal{M}$ and its relation to the quantum cohomology action of this operad upon other total motives. Analogies with homotopy theory, in particular, $A^1$--homotopy formalism, might help to recognize a pattern in algebraic geometry similar to that of iterated loop spaces.

\smallskip

A warning is in order: many meaningful questions cannot be asked and answers cannot be obtained until one extends both parts of the theory, motives and quantum cohomology, from the category $Var_k$ to at least the category of smooth DM--stacks. Some of the arising complications can be avoided if one restricts to the case of genus zero quantum cohomology. We adopt this restriction in this article.

\medskip

{\bf 0.2. Results of this paper.}
This paper is our first installment to the project whose goal is to understand the
Gromov--Witten theory of moduli spaces of curves, preferably on the motivic level,
that is the level of $J$-- and $I$--correspondences (cf. [Beh3] for a nice
and intuitive introduction).

\smallskip

Specifically,  the spaces   $\overline{M}_{0,S}$ (with variable $S$) and their products
are interrelated by a host of natural morphisms expressing
embeddings of boundary strata, forgetting labeled points, relabeling etc:
cf. a systematic description in [BehMa].

\smallskip

Gromov--Witten classes that we study in this paper are certain Chow correspondences 
$$
I(S,\Sigma,\beta ) \in A_*(\overline{M}^{\Sigma}_{0,S}\times \overline{M}_{0,\Sigma})
\eqno(0.1)
$$
where $S$, $\Sigma$ are (disjoint) finite sets of labels and $\beta$ runs over
classes of effective curves in $A_1(\overline{M}_{0,S})$.

\smallskip

Our main motivation is the following vague

\medskip

{\bf 0.2.1. Guess.} {\it Classes (0.1) are ``natural'' in the sense that they 
can be functorially expressed through canonical morphisms in the
category of moduli spaces of labeled trees of various combinatorial types.}

\medskip

This guess is a natural
first step to the understanding the self--referential nature of Gromov--Witten theory in motivic
algebraic
geometry: the fact that components of the basic modular operad
"are" algebras over the same operad (if one takes into account twisting and
grading by the cones of effective curves).

\smallskip

The main result of this paper is an explicit description, in the spirit of our guess,
 of those $I$--correspondences
of $\overline{M}_{0,S}$  that correspond {\it to the classes $\beta$
of boundary curves:} see Theorem 4.6 in section 4.

\smallskip

This answers a question which is quite natural, in particular, because there is a conjecture
that boundary curve classes are generators of the Mori cone: cf. [KeMcK], [FG], [HaT], [CT]
for this and related problems.
 
\medskip

{\bf 0.3. From curves to surfaces and further on?}  One can imaginatively say 
that quantum cohomology of $V$ reveals hidden geometry that can be seen only when
one starts probing $V$ by mapping curves $C$ (``strings'') to $V$. A natural
question arises, how to use, say, surfaces (``membranes'') in place of curves,
and do it in algebraic geometry rather than in symplectic or differential one.

\smallskip

If we expect to discover new universal motivic actions in this way,
we must first contemplate the case $V$= {\it a point} and pose the question:

\smallskip
{\it What are  analogs of moduli spaces $\overline{M}_{g,n}$ (or at least $\overline{M}_{g,0}$) for surfaces
in this context?}

\smallskip

The experience of stringy case indicates that these analogs must be {\it rigid}
objects, as   $\overline{M}_{g,n}$ themselves: see [Hac].

\smallskip

In fact,  moduli spaces are only rarely rigid, but according to a brave
guess of M.~Kapranov, if one starts with an object $X=X(0)$ of dimension $n$, produces
its moduli space of deformations $X(1)$, then moduli space $X(2)$ of deformations
of $X(1)$ etc., then $X(n)$ must be rigid. Quoting [Hac], who summarizes 
philosophy expressed in an unpublished manuscript
by M.~Kapranov,  ``one 
thinks of $X(1)$ as $H^1$ of a sheaf of non--abelian groups on $X(0)$. Indeed, at 
least the tangent space to $X(1)$ at $[X]$ is identified with $H^1(\mathcal{T}_X)$, where $\mathcal{T}_X$ 
is the tangent sheaf, the sheaf of first order infinitesimal automorphisms of 
$X$ . Then one regards $X(m)$ as a kind of non--abelian $H^m$, and the analogy 
with the usual definition of abelian $H^m$ suggests the statement above.''

\smallskip

Extending this idea, one might guess that  an imaginary ``membrane quantum cohomology''
should define motivic actions of rigid (iterated) moduli spaces of surfaces
(endowed with cycles to keep track of incidence conditions) upon certain 
total (ind--)motives. One motivation of this note is to make some propaganda for this idea.

\bigskip

\centerline{\bf 1. Gromov--Witten correspondences}

\medskip

We start with background  terminology and notation.

\medskip
{\bf 1.1. Moduli stacks.} We will consider schemes over a fixed field $k$ of characteristic zero.
\smallskip

Any scheme $W$  ``is the same as" the contravariant functor of its
$T$--points $W(T)=\text{Hom} (T,W)$ with values in $Sets$.
\smallskip

More generally, a stack (of groupoids) $\mathcal{F}$  ``is the same as'' the class of its
$T$--points $\mathcal{F}_T$, where $T$ runs over schemes. The main new element of the situation is that
each $\mathcal{F}_T$ itself and their union $\mathcal{F}$ over ``all'' $T$'s are not simply  sets or  classes, but    categories. Moreover, they form a sheaf on the \'etale (or fppf) site of schemes.

\smallskip

So we will think about
individual objects of $\mathcal{F}_T$ as schematic $T$--points of $\mathcal{F}$, whereas nontrivial morphisms between them
are functors subject to a list of restrictions specific for stacks. Below we recall this list informally.

\smallskip

As in [Ma], V.3, we will imagine objects of $\mathcal{F}_T$ as ``families (of something) over $T$''.
In practical terms, one family is usually given by a diagram of schemes and morphisms,
in which a part of the data remains fixed,
including its  `` base''  $T$, and the rest is subject to a list of explicit restrictions. 

\smallskip

For example, if $\mathcal{F}$ is represented by a scheme $W$, then ``one family over $T$ (of points of $W$)''  is 
a very simple diagram $T\to W$. 

\smallskip

The following requirements must be satisfied.

\medskip

(i) Each object of $\mathcal{F}$ belongs to an $\mathcal{F}_T$ for a unique scheme $T$,
and the map $b:\, \mathcal{F}\to Sch$, sending a family to its base,  is a functor. Groupoid property requires that if $b(f)=id_T$, then 
$f$ is an isomorphism between two respective $T$--points. 

\smallskip

(ii) With respect to morphisms $\varphi :\,T_1\to T_2$,   $\mathcal{F}_T$ must be contravariant:
we must be given  ``base change'' functors $\varphi^*:\,\mathcal{F}_{T_2}\to\mathcal{F}_{T_1}$,
together with functor isomorphisms $(\varphi\circ\psi )^*\to \psi^*\circ \varphi^*$ and associativity diagrams for them.
\smallskip

Moreover, if $F\in Ob\,\mathcal{F}_{T_2}$, then the lifted family $\varphi^*(F)\in Ob\,\mathcal{F}_{T_1}$
must be endowed with a canonical morphism $F\to \varphi^*(F)$ lifting $\varphi$ and
satisfying a set of conditions.

\smallskip

For example, the base change for $T_2\to W$ is  simply the composition $T_1\to T_2\to W$.

\smallskip

(iii) $\mathcal{F}$ is a stack, if each $T$--family is uniquely defined by its restrictions to an \'etale
(or fppf) covering of $T$
and the standard descent data.
The same must be true about morphisms of $T$--families etc. 

\smallskip

(iv) Morphisms of stacks are functors between the respective categories of families, identical on bases
of families.

\smallskip

Thus, an object $F\in Ob\,\mathcal{F}_T$ can also be treated as a
stack, and as such,  it is endowed with a morphism
of stacks $F \Rightarrow \mathcal{F}.$

\medskip

{\bf 1.2. Families of stable maps: preliminaries.}  We will now describe main classes of families
and stacks with which we deal here.

\smallskip

First of all, fix a finite set $\Sigma$,
a genus $g\ge 0$, 
a smooth projective manifold $W$
over $k$ , and an effective class $\beta \in A_1(W)$.

\smallskip
Then  one can define a (proper DM)--stack
$\overline{M}_{g,\Sigma}(W,\beta )$.  

\smallskip

For a $k$--scheme $T$, one object of the groupoid 
$\overline{M}_{g,\Sigma}(W,\beta )(T)$   of  $T$--points
of this stack  consists of a diagram of schemes of the following structure:
\medskip

$$
\xymatrix{\mathcal{C}_T\ar[r]^{f_T} \ar[d]^{h_T} & W\\
 T }
\eqno(1.1)
$$
and a family of sections $x_{j,T}:\, T\to \mathcal{C}_T$, $j\in \Sigma$, 
$h_T\circ x_{j,T}=id_T$.

\smallskip

They must satisfy the following conditions:

\smallskip
(a)  $\mathcal{C}_T\to T$ and $(x_{j,T})$ constitute  a flat  {\it prestable} $T$--family of curves of genus $g.$ 
\smallskip

(b)  $f_T:\, (\mathcal{C}_T; (x_{j ,T}))\to W$, is a {\it stable map} of class $\beta$.

\smallskip
For precise definitions of (pre)stability and class of the map that we use here,  see [BehMa] or [Ma]. 

\smallskip

Given such a diagram with sections, we  call $(W,\beta )$ its {\it target}, $T$ its {\it base}, 
and the whole setup  {\it a $T$--family of stable maps.} Isomorphisms  of families, lifting $id_T$,
must be identical also on $W$. Base changes are defined in a rather evident way.
\smallskip

The stack   $\overline{M}_{g,\Sigma}(W,\beta )$ is defined as the base of
the universal family of this type with given target $(W,\beta )$:
$$
\xymatrix{\overline{C}_{g,\Sigma}(W,\beta )\ar[r]^{\quad f} \ar[d]^{h} & W\\
\overline{M}_{g,\Sigma}(W,\beta )}
\eqno(1.2)
$$
It is endowed with sections $x_j:\,\overline{M}_{g,\Sigma}(W,\beta )  \to \overline{C}_{g,\Sigma}(W,\beta )$,
$j\in \Sigma$.

\smallskip

Naturally,  $\overline{C}_{g,\Sigma}(W,\beta )$ is a stack as well.

\smallskip

If $W$ is a point, $\beta =0$, we routinely omit the target and write simply   $\overline{M}_{g,\Sigma}$,
$\overline{C}_{g,\Sigma}$ etc.
\smallskip
Moreover, (1.2) produces the {\it evaluation/stabilization} diagram
$$
\xymatrix{\overline{M}_{g,\Sigma}(W,\beta )\ar[r]^{\quad st} \ar[d]^{ev} & \overline{M}_{g,\Sigma}\\
W^{\Sigma}}
\eqno(1.3)
$$
Here
$$
ev = (ev_j= f\circ x_j \,|\, j\in \Sigma ):\quad \overline{M}_{g,\Sigma}(W,\beta )\to W^{\Sigma}
\eqno(1.4)
$$
and, in case $2g+|\Sigma |\ge 3$,  the absolute stabilization morphism $st$
discards the map $f$ and stabilizes the remaining prestable family of curves
$$
st:\,\overline{M}_{g,\Sigma}(W,\beta )\to \overline{M}_{g,\Sigma}.
\eqno(1.5)
$$

The {\it virtual fundamental class}, or the $J$--class  $[ \overline{M}_{g,\Sigma}(W,\beta )]^{virt}$,
is a canonical element in the  Chow ring $A_*(\overline{M}_{g,\Sigma}(W,\beta ))$:
$$
J_{g,\Sigma}(W,\beta )\in A_D(\overline{M}_{g,\Sigma}(W,\beta ))\, ,
\eqno(1.6)
$$
where $D$ is the virtual dimension (Chow grading degree) 
$$
(-K_W,\beta ) + |\Sigma |  +(\dim \,W -3)(1-g).
\eqno(1.7)
$$

\smallskip

The respective {\it Gromov--Witten  correspondence}, defined for
$2g+|\Sigma |\ge 3$,    is the proper pushforward
$$
 I_{g,\Sigma}(W,\beta ):=(ev,st)_*( J_{g,\Sigma}(W,\beta )) \in A_D(W^{\Sigma}\times \overline{M}_{g,\Sigma})
 \eqno(1.8)
$$
\smallskip
Understanding  these correspondences is the content of  {\it motivic}
quantum cohomology.

\medskip

{\bf 1.3. Example: $g=0,$ $\beta =0$.} In this case the universal family (1.2) is
$$
\xymatrix{W\times \overline{C}_{0,\Sigma} \ar[r]^{\quad\quad pr_1} \ar[d]^{ id_W\times h} & W\\
W\times \overline{M}_{0,\Sigma} }
\eqno(1.9)
$$
with structure sections $ id _W\times x_j$.

\smallskip

The stabilization morphism is simply the projection
$$
st=pr_2: \  W\times  \overline{M}_{0,\Sigma} \to\overline{M}_{0,\Sigma} .
\eqno(1.10)
$$
\smallskip

The evaluation morphism is the projection followed by the diagonal embedding $\Delta_{\Sigma}$:
$$
ev: \   W\times  \overline{M}_{0,\Sigma} \to W\to W^{\Sigma}.
\eqno(1.11)
$$

We have ([Beh1], p.~606):
$$
J_{0,\Sigma}(W,0)= [\overline{M}_{0,\Sigma}(W,0)]= [W]\otimes [\overline{M}_{0,\Sigma}] .
\eqno(1.12)
$$

\smallskip

The virtual dimension (1.7) is 
$$
|\Sigma| +\dim \, W -3 = \dim \,(W\times \overline{M}_{0,\Sigma}).
$$
Thus, finally, the Gromov--Witten correspondence is the class
$$
 I_{0,\Sigma}(W, 0 ) =
[\Delta_{\Sigma}(W)]\otimes [ \overline{M}_{0,\Sigma}]\in A_*(W^{\Sigma}\times \overline{M}_{0,\Sigma}).
\eqno(1.13)
$$

\medskip
{\bf 1.4. Strategy.} In the remaining sections of this paper, we study the
Gromov--Witten correspondences {\it of genus zero} for $W=\overline{M}_{0,S}$,
$\beta =$ {\it class of a boundary curve}  in  $\overline{M}_{0,S}$ (cf. below).
This is the first step of a much more ambitious program in which all components of the
stable family diagrams are allowed to be stacks, and and in which we take for
targets  the stacks $\overline{M}_{g,S}$ and arbitrary $\beta$.

\smallskip

Our modest goal here allows us to basically restrict ourselves to the case of schemes,
whose geometry is already well known.  However, some  intermediate constructions
require the use of stacks.

\smallskip

In particular , we need to understand  the relevant $J$--classes and the diagrams
$$
ev:\,\overline{M}_{0,\Sigma}(\overline{M}_{0,S},\beta  )\to \overline{M}^{\Sigma}_{0,S},\quad
st:\,\overline{M}_{0,\Sigma}(\overline{M}_{0,S},\beta)\to \overline{M}_{0,\Sigma}.
\eqno(1.14)
$$

\smallskip

We also want to be able to trace various functorialities, in particular, in {\it both} $S$ and $\Sigma$.
However, this may result in a rather clumsy notation.

\smallskip

In order to postpone its introduction, in the remaining parts of this
section we describe a somewhat more general  situation. Afterwards we will show that
our main problem is contained in it.

\medskip

{\bf 1.5. Setup, part I.} Consider a morphism of smooth irreducible projective manifolds  $b:\, E\to W$. Let $\beta_{E}$ be
an effective curve class on $E$, and $\beta := b_*(\beta_E)$ its pushforward to $W$. Any stable map
$\mathcal{C}_T/T\to E, (x_j:T\to \mathcal{C}_T\,|\,j\in \Sigma )$,  of class $\beta_E$, induces, after composition with $b$ and stabilization, a stable map with target $(W,\beta )$. Thus, we get a map
$$
\widetilde{b}:\ \overline{M}_{0,\Sigma}(E,\beta_E)\to   \overline{M}_{0,\Sigma}(W,\beta)
$$
that clearly fits into the commutative diagram
 $$
\xymatrix{\overline{M}_{0,\Sigma}(E,\beta_E ) \ar[r]^{\widetilde{b}}     \ar[d]^{(ev_E,st_E)}& 
 \overline{M}_{0,\Sigma}(W,\beta )\ar[d]^{(ev_W,st_W)} \\
E^{\Sigma} \times \overline{M}_{0, \Sigma} \ar[r]^{b^{\Sigma}\times id}& W^{\Sigma}\times \overline{M}_{0,\Sigma}}
\eqno(1.15)
$$

If $|\Sigma |\le 2$, the space  $\overline{M}_{0,\Sigma}$ is not a DM--stack;
discarding it and stabilization morphisms in (1.15), we still get a commutative diagram.
Whenever  $\overline{M}_{0,\Sigma}$ appears, we assume that $|\Sigma |\ge 3$.
\medskip

{\bf 1.5.1. Proposition.} {\it (i) Assume that
$$
J_{0,\Sigma}(W,\beta )=\widetilde{b}_{*}(J_{0,\Sigma}(E,\beta_E)).
\eqno(1.16)
$$
Then
$$
I_{0,\Sigma}(W,\beta )=(b^{\Sigma}\times id)_{*}(I_{0,\Sigma}(E,\beta_E)).
\eqno(1.17)
$$

\smallskip

(ii) Let $\gamma_j \in H^*(W),$ $j\in \Sigma$, be a family
of cohomology classes marked by $\Sigma$. Then from (1.16) it follows that
$$
pr^*_W(\otimes_{j\in \Sigma} \gamma_j)\cap I_{0,\Sigma}(W,\beta ) =
$$
$$
=(b^{\Sigma}\times id)_* [pr^*_E(\otimes_{j\in \Sigma} b^*(\gamma_j))\cap I_{0,\Sigma}(E,\beta_E)] .
 \eqno(1.18)
$$
Here we denote by $pr_W:\,W^{\Sigma} \times \overline{M}_{0, \Sigma}\to W^{\Sigma}$
and  $pr_E:\,E^{\Sigma} \times \overline{M}_{0, \Sigma}\to E^{\Sigma}$
the respective projection morphisms, and $H^*$ can be any standard cohomology theory.
}

\medskip

{\bf Proof.} (i) This follows directly from (1.16) and commutativity of  (1.15).

\medskip

(ii) We have, using the projection formula
$$
(b^{\Sigma}\times id)_* [pr^*_E(\otimes_{j\in \Sigma} b^*(\gamma_j))\cap I_{0,\Sigma}(E,\beta_E)]  =
$$
$$
=(b^{\Sigma}\times id)_*[ (b^{\Sigma}\times id)^* \circ pr^*_W(\otimes_{j\in \Sigma} \gamma_j)\cap I_{0,\Sigma}(E,\beta_E)]=
$$
$$
= pr^*_W(\otimes_{j\in \Sigma} \gamma_j)\cap (b^{\Sigma}\times id)_*(I_{0,\Sigma}(E,\beta_E))
 $$
The last expression coincides with l.h.s. of (1.18) in view of (1.17). This completes the proof.

\medskip

{\bf 1.5.2. Remark.} In our applications to the case $W=\overline{M}_{0,S}$
(cf. section 4), $E$ will be a boundary stratum
containing the boundary curve representing $\beta$,  and the virtual fundamental classes $J_{0,\Sigma}$
will coincide with the usual fundamental classes since the relevant deformation
problem will be unobstructed. Moreover, $E$ will have a very special additional structure.
We will axiomatize the relevant  geometry in the next subsections.

\medskip

{\bf 1.6. Setup, part II.}  Keeping notation  of  1.5, we  make the following additional
assumptions:

\medskip

{\it (a) $E$ is explicitly represented as $E=B\times C$ where $C$ is isomorphic to $\textbf{P}^1$.
This identification, including the projections $p=pr_B:\, E\to B$ and $pr_C:\, E\to C$,
constitutes  a part of structure.

\medskip

(b) $\beta_E$ is the (numerical) class of any fiber of $p$.

\medskip

(c) The deformation problem for any fiber $C_0$ of $p$ embedded via $b_0$ in $W$
is {\it trivially unobstructed} in the sense of  [Beh3]:
$$
H^1(C_0, b_0^*(\mathcal{T}_W))=0 \, .
\eqno(1.19)
$$

(d) The map $\widetilde{b}$ in (1.15) is an isomorphism. }

\medskip

These assumptions are quite strong.  In particular, from (b) -- (d)
it follows that (1.16) holds since the relevant virtual fundamental classes coincide with the ordinary
ones.
Thus, we can complete the explicit computation of
$I_{0,\Sigma} (W,\beta )$ starting with the right hand side of (1.17).
We will do it in the remaining part of the section.

\medskip

First of all, we have
$$
pr_{B*}(\beta_E) =0,\quad pr_{C*}(\beta_E)= \textbf{1} \,
$$
where    $\textbf{1}$ is the fundamental class  $[C]$ in the Chow ring of $C$.
\smallskip

Thus, the two projections induce the map
$$
(\widetilde{pr}_B,\widetilde{pr}_C):
\overline{M}_{0,\Sigma}(E,\beta_E)\to  \overline{M}_{0,\Sigma}(B,0)\times 
\overline{M}_{0,\Sigma}(C,\textbf{1})\, .
$$
Stabilization maps embed this morphism into  the commutative diagram
$$
\xymatrix{\overline{M}_{0,\Sigma}(E,\beta_E)\ar[r] \ar[d]_{st_E} & \overline{M}_{0,\Sigma}(B,0)\times 
\overline{M}_{0,\Sigma}(C,\textbf{1}) \ar[d]_{st_B\times st_C} \\
\overline{M}_{0,\Sigma}\ar[r]^{\Delta_{\overline{M}_{0,\Sigma}}} & \overline{M}_{0,\Sigma}\times \overline{M}_{0,\Sigma}}
\eqno(1.20)
$$
where the lower line is the diagonal embedding (cf. [Beh2], Proposition 5).

\smallskip

Similarly, evaluation maps   embed this morphism into  the commutative diagram
$$
\xymatrix{\overline{M}_{0,\Sigma}(E,\beta_E)\ar[r] \ar[d]_{ev_E} & \overline{M}_{0,\Sigma}(B,0)\times 
\overline{M}_{0,\Sigma}(C,\textbf{1}) \ar[d]_{ev_B\times ev_C} \\
E^{\Sigma}\ar[r]^{s} & B^{\Sigma}\times C^{\Sigma}}
\eqno(1.21)
$$
where the lower line is now the evident permutation isomorphism induced by
$E=B\times C$.

\smallskip

Combining these two diagrams, we get
$$
\xymatrix{\overline{M}_{0,\Sigma}(E,\beta_E)\ar[r] \ar[d]_{(ev_E,st_E)} & \overline{M}_{0,\Sigma}(B,0)\times 
\overline{M}_{0,\Sigma}(C,\textbf{1}) \ar[d]^{(ev_B,st_B)\times (ev_C,st_C)} \\
E^{\Sigma}\times\overline{M}_{0,\Sigma}\ar[r]^{\widetilde{\Delta}\quad\quad\quad} & B^{\Sigma}\times \overline{M}_{0,\Sigma}
\times C^{\Sigma}\times \overline{M}_{0,\Sigma}}
\eqno(1.22)
$$
Here the lower line is an obvious composition of permutations and the diagonal embedding of 
$\overline{M}_{0,\Sigma}$.

\medskip

From (1.22) and [Beh2]  it follows that
$$
 I_{0,\Sigma}(E, \beta_E ) = \widetilde{\Delta}^! (I_{0,\Sigma}(B,0)\otimes I_{0,\Sigma}(C,\textbf{1}))\, .
 \eqno(1.23)
$$
(Notice that for $x\in A_*(X), y\in A_*(Y)$ we often denote simply by
$x\otimes y\in A_*(X\times Y)$ the image of  $x\otimes y\in A_*(X)\otimes A_*(Y)$
wrt the canonical map $ A_*(X)\otimes A_*(Y)\to A_*(X\times Y)$).

\smallskip
Furthermore, according to (1.13), 
$$
 I_{0,\Sigma}(B, 0 ) =
[\Delta_{\Sigma}(B)\times  \overline{M}_{0,\Sigma}]\in A_*(B^{\Sigma}\times \overline{M}_{0,\Sigma}).
\eqno(1.24)
$$

Finally, the space $\overline{M}_{0,\Sigma}(C,\textbf{1})$ and the class 
$I_{0,\Sigma}(C,\textbf{1})$ can be described as follows.

\smallskip
Recall a construction from [FuMPh].
Let $V$ be a smooth complete algebraic manifold. For a finite set $\Sigma$,
let $V^{\Sigma}$ be the direct product of a family of $V$'s labeled by elements of $\Sigma$.
Denote by $\widetilde{V}^{\Sigma}$ the blow up
of the (small) diagonal in $V^{\Sigma}$. Finally, define $V^{\Sigma ,0}$
as the complement to all partial diagonals in   $V^{\Sigma}$. 

\smallskip

The Fulton--MacPherson's {\it configuration space $V\langle \Sigma \rangle$}
(for curves it was earlier introduced by Beilinson and Ginzburg) is the closure
of  $V^{\Sigma ,0}$ naturally embedded in the product
 $$
 V^{\Sigma}\times\prod_{\Sigma^{\prime}\subset \Sigma , |\Sigma^{\prime}|\ge 2}
\widetilde{V}^{\Sigma^{\prime}} .
$$
 In [FuPa], it was shown that  $\overline{M}_{0,\Sigma}(C,\textbf{1})$ can be identified
 with $C\langle \Sigma\rangle$ in such a way that the  birational morphism
 $ev_C$  becomes the tautological open embedding  when restricted
 to $C^{\Sigma ,0}$.  
 \smallskip
 
 Therefore, denoting by $D_{\Sigma}\subset C^ \Sigma \times \overline{M}_{0,\Sigma}$
 the closure of the graph of the  canonical surjective map
$C^{\Sigma ,0} \to   M_{0,\Sigma}$, we get 
$$
I_{0,\Sigma}(C,\textbf{1}) = [D_{\Sigma}] \, .
\eqno(1.25)
$$

Now we can state the main result of this section:

\medskip

{\bf 1.6.1. Proposition.} {\it Assuming 1.6 (a)--(d), we have
$$
 I_{0,\Sigma}(E, \beta_E ) = \widetilde{\Delta}^! ([\Delta_{\Sigma}(B)\times  \overline{M}_{0,\Sigma}\times 
 D_{\Sigma}])\, .
 \eqno(1.26)
$$
and
$$
I_{0,\Sigma}(W,\beta )=(b^{\Sigma}\times id)_{*}\circ
\widetilde{\Delta}^! ([\Delta_{\Sigma}(B)\times  \overline{M}_{0,\Sigma}\times 
 D_{\Sigma}])\, .
\eqno(1.27)
$$

}
This is a straightforward  combination of  (1.23) -- (1.25) and (1.17).

\bigskip

\centerline{\bf 2. Target space $\overline{M}_{0,4}$}

\medskip

{\bf 2.1. Notation.} Stressing functoriality wrt labeling sets, and having in mind 
further developments, we denote in this section by $S$  a set {\it of cardinality 4},
with a marked point $\bullet$. We put $S=P\sqcup \{\bullet \}$.
Thus, we are considering the moduli space
$\overline{M}_{0,P\sqcup \{\bullet \}}$.  It is a projective line
endowed with  three pairwise distinct  points $D_{\sigma}$ labeled by  unordered partitions
$\sigma :\,P\sqcup \{\bullet\}=S_1\sqcup S_2$, $|S_i|=2.$ They are exactly those points
over which the universal stable curve $\overline{C}_{0,P\sqcup\{\bullet \} }$ splits into two components,
and labeled points are redistributed among them according to $\sigma$. 
Now, the set of such partitions is naturally bijective to $P$: $j \in P$ corresponds
to the partition $\{\bullet, j \}\sqcup (P\setminus\{\bullet ,j \})$. Hence finally
$\overline{M}_{0,P\sqcup \{\bullet \}}$  is a projective line $\textbf{P}^1$ stabilized by three points labeled by $P$.
This identification is functorial wrt pointed bijections of $S$.

\smallskip

The only boundary class  of curves in $A_1( \overline{M}_{0,P\sqcup \{\bullet\}})$ is the fundamental class
$\beta =\textbf{1}: =  [\overline{M}_{0,P\sqcup \{\bullet\}}]$. 
We have already invoked the description
of the universal family of stable maps with this target and the relevant $I$--class at the end of 
1.6, see (1.25). But now, for the sake of a future generalization, we will use a slightly different 
family and an alternative
description  of  $I_{0,\Sigma}\in A_*( ( \textbf{P}^1)^{\Sigma}\times \overline{M}_{0,\Sigma})$
that will better fit the passage to target spaces $\overline{M}_{0,S}$ with
$|S|>4.$

 \medskip

{\bf 2.2. An alternative family.} Consider
the moduli space $\overline{M}_{0,\Sigma\sqcup P\sqcup \{\bullet\}}$.

\smallskip
Recall that for any finite set $R$ and its subset $Q\subset R$ with complement
of cardinality $\ge 3$,
the space $\overline{M}_{0,R}$ is the source of the
standard {\it forgetful morphism} $\psi_Q :\, \overline{M}_{0,R}\to \overline{M}_{0,R\setminus Q}$: 
 ``forget the subset of sections labeled by $Q$ and stabilize".
\smallskip

Thus we get the diagram
$$
\xymatrix{ \overline{M}_{0,\Sigma\sqcup P\sqcup \{\bullet\}}\ar[r]^{\psi_{\Sigma}} \ar[d]^{\psi_{\{\bullet\}}} 
& \overline{M}_{0,P\sqcup \{\bullet\}}\\
 \overline{M}_{0,\Sigma\sqcup P}  }
\eqno(2.1)
$$

Another standard morphism identifies   the vertical arrow in (2.1)
with the projection of the universal $(\Sigma\sqcup P)$--labeled curve 
$\overline{C}_{0,\Sigma\sqcup P}$ to its base: cf. e.g. [Ma], Ch. V, Theorem 4.5.

\smallskip

From the explicit form of this identification, one easily sees that the image 
in  $\overline{M}_{0,\Sigma\sqcup P\sqcup \{\bullet\}}$ of the section
$x_j:\,\overline{M}_{0,\Sigma\sqcup P} \to  \overline{C}_{0,\Sigma\sqcup P}$ for $j\in \Sigma\sqcup P$
is precisely the boundary divisor of  $\overline{M}_{0,\Sigma\sqcup P\sqcup \{\bullet \}}$
corresponding to the stable 2--partition 
$$
\Sigma\sqcup P\sqcup \{\bullet\}= \{\bullet ,j\}\sqcup ((\Sigma\sqcup P)\setminus\{j\})
\eqno(2.2)
$$
Here we will denote this divisor by $D_j$.

\smallskip
Consider now (2.1) as the family of maps of class $\textbf{1}$, in which {\it only the sections
$x_j$ for $j\in \Sigma$ are counted as structure sections,} whereas those labeled by $P$
are discarded. Then the family will not be stable anymore: if an irreducible component
of fiber curve contains only three special points and one of them corresponds to the section
labeled by an element of $P$, then this component will be contracted by $\psi_{\Sigma}$.
We can stabilize this new family and get a diagram
$$
\xymatrix{\overline{C}\ar[r]^{\overline{\psi}_{\Sigma}\quad} \ar[d] &\overline{M}_{0,P\sqcup \{\bullet\}}\\
 \overline{M}_{0,\Sigma\sqcup P}  }
\eqno(2.3)
$$
endowed additionally with sections labeled by $\Sigma$ and the stabilizing morphism
$$
\chi:\ \overline{M}_{0,\Sigma\sqcup P\sqcup \{\bullet\}}\to \overline{C}, \quad 
{\psi}_{\Sigma} =\overline{\psi}_{\Sigma}\circ \chi.
\eqno(2.4)
$$
\smallskip

For each $j\in \Sigma$, denote by  $\xi_j :\, \overline{M}_{0,\Sigma\sqcup P} \to 
\overline{M}_{0,\Sigma\sqcup P\sqcup \{\bullet\}}$
the section of $\psi_{\{\bullet\}}$ identifying  $\overline{M}_{0,\Sigma\sqcup P}$ with 
$D_j\subset  \overline{M}_{0,\Sigma\sqcup P\sqcup \{\bullet\}}$ from (2.2).
Consider the map
$$
\overline{ev} := (\psi_{\Sigma}\circ \xi_j\,|\,j\in\Sigma ):\ \overline{M}_{0,\Sigma\sqcup P}
\to   (\overline{M}_{0,P\sqcup \{\bullet\}})^{\Sigma}.
\eqno(2.5)
$$

\smallskip

The stable family (2.3) may be obtained by a base change from the universal family of stable maps of class 
$\beta$. Let
$$
\mu :\ \overline{M}_{0,\Sigma\sqcup P} \to \overline{M}_{0,P\sqcup \{\bullet\}}\langle \Sigma\rangle 
\eqno(2.6)
$$
be the respective  morphism of bases.
\medskip

Dimensions of the two smooth irreducible schemes in (2.6) coincide.
It is not difficult to see that morphism $\mu$ is  birational and hence surjective. In fact, consider a 
generic fiber of  $\overline{C}_{0,\Sigma\sqcup P}$. It is simply $\textbf{P}^1$
with pairwise distinct  $(\Sigma\sqcup P(\Pi ))$--labeled points. When we discard $\Sigma$--labeled
ones, we land in  $\textbf{P}^1$ endowed with three points labeled by $P(\Pi )$; inverse images of them 
are just missing section that we discarded when constructing (2.3) from (2.1);
so in fact at a generic point we neither loose, nor gain any information passing from (2.1)
to (2.3).

\smallskip

We can now prove the main result of this section.

\medskip

{\bf  2.3. Proposition.} {\it We have for $|\Sigma|\ge 3$:
$$
I_{0,\Sigma}(\overline{M}_{0,P\sqcup \{\bullet \}},\textbf{1} )= (\overline{ev}, \psi_P)_*
([\overline{M}_{0,\Sigma\sqcup P}])\in A_{|\Sigma |}( (\overline{M}_{0,P\sqcup \{\bullet \}})^{\Sigma}
\times \overline{M}_{0,\Sigma}    ).
\eqno(2.7)
$$
}
\smallskip

{\bf Proof.}  Since $\mu$ is birational and surjective, we can identify
the relevant $J$--class with
$$
\mu_*([ \overline{M}_{0,\Sigma\sqcup P}])= [ \overline{M}_{0,P\sqcup \{\bullet \}}\langle \Sigma\rangle ].
$$
In order to prove (2.7), it remains to check that
$$
ev\circ \mu= \overline{ev},\quad  st\circ \mu =\psi_P.
\eqno(2.8)
$$
Both facts follow from the discussion in 2.2 above.

\bigskip

\centerline{\bf  3.  Boundary curve classes in $\overline{M}_{0,S}$}

\medskip

In this section, after recalling some basic facts about boundary of $\overline{M}_{0,S}$ following [Ma]
and [BehMa], we  
summarize relevant parts of [KeMcK] and fix  our notation.

\medskip

{\bf 3.1.  Boundary strata of   $\overline{M}_{0,S}$.}  
The main combinatorial invariant of an $S$--pointed stable curve $C$ is its {\it dual graph} $\tau =\tau_C$.
Its set of vertices $V_{\tau}$ is (bijective to) the set of irreducible components of $C$.
Each vertex $v$ is a boundary point of the set of {\it flags}  $f\in F_{\tau}(v)$ which is (bijective to)
the set consisting of singular points and $S$--labeled points on this component. We put 
$F_{\tau}=\cup_{v\in V_{\tau}}F_{\tau}(v).$
If two components of $C$ intersect, the respective two vertices carry two flags
that are grafted to form an {\it edge} $e$ connecting the respective vertices; the set of edges is denoted $E_{\tau}$.
The flags that are not pairwise grafted are called {\it tails.} They form a set $T_{\tau}$
which is naturally bijective to the set of $S$--labeled points and therefore itself is
labeled by $S$.   Stable curves of genus zero correspond to trees
$\tau$ whose each vertex carries at least three flags.

\smallskip

The space $\overline{M}_{0,S}$ is a disjoint union of locally closed strata $M_{\tau}$
indexed by stable $S$--labeled trees. Each such stratum $M_{\tau}$ represents the functor
of families consisting of curves of combinatorial type $\tau$. In particular, the open stratum
$M_{0,S}$ classifies irreducible smooth curves with pairwise distinct $S$--labeled points.
Its graph is a star:  tree with one vertex, to which all tails are attached, and having no edges.

\smallskip

Generally, a stratum $M_{\tau}$ lies in the closure $\overline{M}_{\sigma}$ of $M_{\sigma}$,
iff $\sigma$ can be obtained from $\tau$ by contracting a subset of edges.
Closed strata $ \overline{M}_{\sigma}$ corresponding to trees with nonempty set of edges,
are called {\it boundary} ones. The number of edges is the codimension of the stratum.

\medskip

{\bf 3.1.1. Boundary divisors and $A^1(\overline{M}_{0,S})$.} The classes of boundary divisors 
generate the whole Chow ring, but are not linearly independent.
The following useful basis is constructed in [FG].

\smallskip

For  $s\in S$, let $\mathcal{L}_s$ be the line bundle on $\overline{M}_{0,S}$ whose
fiber over a stable curve $(C,(x_t\,|\,t\in S))$ is $T^*_{x_s}C$. Put  
$\psi_s:=c_1(\mathcal{L}_s  )$

\medskip

{\bf 3.1.2. Proposition.}  {\it The classes of boundary divisors 
$D_S$ with
$|S_1|, |S_2| \ge 3$ and classes $\psi _s$,  $s\in S$,
constitute a basis of the group $A^1(\overline{M}_{0,S})$. 

\smallskip

The rank of this group is $2^{n-1} -\dfrac{n(n-1)}{2} -1.$}
 
\medskip

This is  Lemma 2 in [FG]. An expression of $\psi_s$ through boundary classes
is given in Lemma 1 of [FG].

\medskip
Below we give some details  on one--dimensional strata.

\medskip

{\bf 3.2. Boundary curves and $A_1(\overline{M}_{0,S})$: preparatory combinatorics.} We start with the following
combinatorial construction.  

\smallskip

We will use here the term {\it an unordered partition} of a set $S$ as synonymous to
{\it an equivalence relation} on $S$. A component of a partition is the same as an equivalence class
of the respective relation;
in particular, all components are non--empty.

\smallskip

Call an unordered partition $\Pi$ of $S$  {\it distinguished}, if it consists of
precisely {\it four  components}. Denote by the $S(\Pi )$ the set
of the components, that is, the quotient of $S$ wrt the respective equivalence relation.

\smallskip

Distinguished partitions are in a natural bijection with isomorphism classes of {\it distinguished stable
$S$--labeled trees $\pi$}.  By definition, such a tree is endowed with one {\it distinguished vertex
$v_0$}, the set of flags at this vertex $F_{\pi}(v_0)$ being (labeled by) elements of $S(\Pi )$.
Clearly this vertex is of multiplicity four.
The flags labeled by one--element components $\{s\}$ of $\Pi$ are tails, carrying the respective 
labels $s\in S$. The remaining flags are halves of edges; the second vertex of an edge, whose one
half is labeled by a component $S_i$ carries tails labeled by elements of $S_i$.

\smallskip

We will routinely identify $F_{\pi}(v_0)$ with  $S(\Pi )$.

\medskip

{\bf 3.2.1. Definition.} {\it (i) Given a distinguished partition $\Pi$,  denote by
$P=P(\Pi )$  the set of those stable 2--partitions $\sigma$ of $S$, each component of which
is a union of two different components of $\Pi $. For $|S|\ge 4$ we have $|P(\Pi )|=3.$

\smallskip

(ii)  $N = N(\Pi )$ is the set of those stable 2--partitions of $S$ whose one component
coincides with one component of $\Pi$. We have for $|S|\ge 5$:
$1\le |N(\Pi )|\le 4.$}

\medskip
{\bf 3.2.2. Lemma.} {\it $\Pi$ can be uniquely reconstructed from $P(\Pi)$;
hence $P(\Pi )$ uniquely determines $N(\Pi )$ as well.}

\smallskip

{\bf Proof.}  In fact, if $\Pi = (S_1, S_2, S_3, S_4)$ (numeration arbitrary), then by definition
$P(\Pi )$ must consist of partitions
$$
\sigma_1=(S_1\cup S_2, S_3\cup S_4),\  
\sigma_2=(S_1\cup S_3, S_2\cup S_4),\  
\sigma_3=(S_1\cup S_4, S_2\cup S_3)
$$
Hence conversely, knowing $P(\Pi )$, we can unambiguously reconstruct
$\Pi $: its components are exactly non--empty pairwise intersections
of components of different $\sigma_i\in P(\Pi ).$

\medskip

{\bf 3.3. Boundary curves and $A_1(\overline{M}_{0,S})$: geometry.}  Each 
distinguished partition $\Pi$
of $S$ determines the following boundary stratum of  $\overline{M}_{0,S}$:
$$
b_{\Pi}:\quad \overline{M}_{\Pi} := \cap_{\sigma\in N(\Pi)} D_{\sigma}\hookrightarrow \overline{M}_{0,S} \, .
\eqno(3.1)
$$
Equivalently, $\overline{M}_{\Pi}$ is the stratum, corresponding to the special tree $\pi$  associated to $\Pi$.
In other words, now all components of $\Pi$ are  indexed by the flags $f\in F_{\pi}(v_0)$ at the special vertex $v_0$, whereas components of cardinality $\ge 2$ are also naturally indexed 
by the remaining vertices of $\pi$: 
$$
\overline{M}_{\Pi} =  \overline{M}_{0,F_{\pi}(v_0)} \times \prod_{v\ne v_0} \overline{M}_{0,F_{\pi}(v)}.
\eqno(3.2)
$$
Here the equality sign refers to the canonical isomorphism that is defined for any stable marked tree:
it produces from such a tree the product of moduli spaces corresponding to the stars of all vertices.
\smallskip

The information about edges determines the  embedding morphism (3.1) of such a product
as a boundary stratum.
On the level of universal curves, it is defined by merging 
the pairs of sections labeled by halves of an edge.

\smallskip

Codimension of $\overline{M}_{\Pi}$ is $|N(P)|$, and $1\le |N(\Pi)| \le 4.$  Since $|F_{\pi}(v_0)|=4$, 
the moduli space
$\overline{M}_{0,F_{\pi}(v_0)}$  is $\textbf{P}^1$ with three points naturally labeled
by the set of  stable partitions of  $F_{\pi}(v_0)$ which in turn is canonically bijective to $P(\Pi)$,
cf. 2.1.

\medskip

Hence the representation (3.2) allows us to define the projection map
$$
p=p_{\Pi}: \, \overline{M}_{\Pi} \to B_{\Pi} := \prod_{v\ne v_0} \overline{M}_{0,F_{\pi}(v)}
\eqno(3.3)
$$
having three canonical disjoint sections canonically labeled by   $P(\Pi)$.

\smallskip

Clearly, all fibers of $p_\Pi $ are rationally equivalent so that they define
a class $\beta = \beta (\Pi )\in A_1(\overline{M}_{0,S})$.

\medskip

{\bf 3.3.1. Lemma} ([KeMcK]). {\it (i) For $n:=|S| \ge 4$, each boundary curve (one--dimensional boundary stratum) $C_{\tau}$ 
is a fiber of  one of the projections (3.3).

\smallskip

(ii)  $[C_{\tau_1}]=[C_{\tau_2}] \in A_1(\overline{M}_{0,S})$ iff these curves  are fibers
of one and the same projection (3.3).}

\medskip
We reproduce the proof for further use.

\smallskip

{\bf Proof.}  (i) Since $C_{\tau}$ is a curve, the $S$--labeled stable tree
$\tau$ is a tree with $|E_{\tau}|=n-4$ and hence
$|V_{\tau}|=n-3.$ Since the tree is stable, all except one of its vertices must have multiplicity
 3. The exceptional vertex  denoted $v_0=v_0(\tau )$ has multiplicity 4.

\smallskip

If we delete from the geometric tree $\tau$ the vertex $v_0$, it will break into 4 connected components.
Thus, the set $S$ of labels of tails will be broken into 4 non--empty subsets. Among them
there are $|T_{\tau}(v_0)|$ one--element sets
(labels of tails adjacent to $v_0$), and $|E_{\tau}(v_0)|$ sets of cardinality $\ge 2$:
each part consist of labels of those tails that can be reached
from the critical vertex by a path (without backtracks) starting with the respective flag.
We will denote this partition $\Pi (\tau )$. Hence if we contract  all edges of $\tau$
excepting those that are attached to $v_0$, we will get the distinguished tree associated
with a distinguished partition $\Pi = \Pi (\tau ).$ It determines the required projection.

\smallskip

(ii) Now consider two sets of stable 2--partitions of $S$ produced from $\Pi=\Pi(\tau )$
as in the Definition 3.2.1, and denote them respectively $P(\tau )$
and $N(\tau )$.

\smallskip

First of all, we will check that 
$$
(D_{\sigma}, C_{\tau})=1,\  \text{if}\   \sigma \in P(\tau ),
$$
$$
(D_{\sigma}, C_{\tau})=-1,\  \text{if}\   \sigma \in N(\tau ),
\eqno(3.4)
$$
$$
(D_{\sigma}, C_{\tau})=0\quad  \text{otherwise.}
$$
Now we will use formulas and facts proved in [Ma], III.3 and [KoMaKa], Appendix.
In particular, we use the notion of {\it good monomials}, elements of the commutative
polynomial ring freely generated by symbols
$m(\sigma )$ where $\sigma$ runs over stable 2--partitions of $S$.
These monomials form a family indexed by stable $S$--labeled trees $\tau$:
$m(\tau ):= \prod_{e\in E_{\tau}} m(\sigma_e)$ where $\sigma_e$ is the 2--partition of $S$
obtained by cutting $e$.

\smallskip

Assume first that $m (\sigma )m(\tau )$ is a good monomial so that
$(D_{\sigma},C_{\tau})=1.$ Then it is of the form $m(\rho )$
where $\rho$ is a stable $S$--labeled tree with all vertices of multiplicity 3 and
an edge $e$ such that  $m(\sigma ) = m(\rho_e)$.  This edge is unambiguously characterized
by the fact that after collapsing $e$ in $\rho$ to one vertex, we get  the labeled tree (canonically
isomorphic to)  $\tau$. But the vertex to which $e$ collapses must
then have multiplicity larger than 3. It follows that $e$ must collapse precisely to the exceptional vertex
$v_0$ of $\tau$.  Conversely, the set of ways of putting $e$ back is clearly in a  bijection with $P(\tau )$:
the 4 flags adjacent to $v_0$ must be distributed in two groups, 2 flags in each, that will
be adjacent to two ends of $e$.

\smallskip

Assume now that $m(\sigma )$ divides $m(\tau )$. Using Proposition 1.7.1 of [KoMaKa],
one sees that $m(\sigma )m(\tau )$ represents zero in the Chow ring (and so $(D_{\sigma},C_{\tau})=0$)
unless $\sigma = \tau_e$ where $e$ is an edge adjacent to $v_0$. In this latter case Kaufmann's formula (1.9)
from [KoMaKa] implies   $(D_{\sigma},C_{\tau})= -1$. The set of such $\sigma$'s is in a bijection with $N(\tau )$.

\smallskip

Finally, for any other stable  2--partition $\sigma$ there exists an $e\in E_{\tau}$ such that
we have $a(\sigma , \tau _e)=4$ in the sense of [Ma1], III.3.4.1. In this case,  $(D_{\sigma},C_{\tau})=0$
in view of  [Ma], III.3.4.2.

\smallskip

Now, we have   $[C_{\tau_1}]= [C_{\tau_2}]$   iff 
 $(D_{\sigma}, C_{\tau_1})=(D_{\sigma}, C_{\tau_2})$ for all stable 2--partitions $\sigma$,
 because boundary divisors generate $A^1$. In view of (3.4),
 this latter condition means precisely  that 
 $$
 P(\tau_1) = P(\tau_2),\quad N(\tau_1) = N(\tau_2).
 $$
But lemma 3.2.2 shows that in this case $\Pi (\tau_1)=\Pi (\tau_2)$. This completes the proof.

\medskip

{\bf 3.4. Proposition.} {\it  Denote the canonical class of  $\overline{M}_{0,S}$ by $K_S$.
Using notation of 3.3, we have
$$
(-K_S,\beta (\Pi ))= 2-|N(\Pi )|.
\eqno(3.5)
$$
}
\medskip

{\bf Proof.} For $2\le j \le [n/2]$, denote by $B_j$ the sum of all divisors $D_{\sigma}$
such that one part of the partition $\sigma$ is of cardinality $j$, and by $B$ the sum of all boundary divisors.  We have
$$
-K_S = 2B-\sum_{j=2}^{[n/2]} \frac{j(n-j)}{n-1}B_j
\eqno(3.6)
$$
(cf.  [KeMcK], [FG], and references therein).
\smallskip

For a stable 2--partition $\sigma =(S_1, S_2)$ of $S$, put
$c(\sigma ) := |S_1| |S_2|.$ Then, combining (3.4) and (3.6), we get:
$$
(-K_S, \beta (\Pi ))=2(|P(\tau )|-|N(\tau )|) - \sum_{\sigma\in P(\tau )}  \frac{c(\sigma )}{n-1} +
\sum_{\sigma\in N(\tau )}  \frac{c(\sigma )}{n-1} .
\eqno(3.7)
$$
The most straightforward way to pass from (3.7) to (3.5) is to consider the four cases
$|N(\Pi )|=1,2,3,4$ separately. Here is the calculation for  $|N(\Pi )|=3$;
it demonstrates the  typical cancellation pattern. We leave the remaining cases to the reader.

\smallskip
We have  $2(|P(\Pi )|-|N(\Pi )|) = 0.$ Let $(1,a,b,c)$ be the cardinalities of
the components of $\Pi$, where $a,b,c\ge 2$, $a+b+c=n-1$.
Then $P(\Pi )$ consists of three partitions, of the following cardinalities respectively
$$
(a+1, b+c),\ (b+1,a+c),\  (c+1, a+b).
$$
Hence
$$
\sum_{\sigma\in P(\Pi )}  c(\sigma ) = 2(ab+ac+bc)+2(a+b+c).
$$
Similarly, partitions in $N(\Pi )$ produce the list
$$
(a,1+ b+c),\ (b,1+a+c),\  (c,1+ a+b)
$$
so that
$$
\sum_{\sigma\in N(\Pi )}  c(\sigma ) = 2(ab+ac+bc)+(a+b+c).
$$
Combining all together, we get $(-K_S,\beta (\Pi ))= -1=2-|N(\Pi )|.$

\medskip

{\bf 3.5. Proposition.} {\it Each  class of a boundary curve  $\beta$ is indecomposable
in the cone of effective curves.}

\medskip

{\bf Proof.} This follows from (3.5) and [KeMcK], Lemma 3.6: $(K_S+B,\beta (\Pi ))=1$, and
the divisor $K_S+B$ is ample.

\medskip

{\bf 3.6. Examples: $\overline{M}_{0,4}$ and $\overline{M}_{0,5}$.}  If $|S|=4$, there is one distinguished partition $\Pi$,
with all components of cardinality 1. The respective ``boundary'' curve is in fact
the total space $\overline{M}_{0,S}.$

\smallskip

 If $|S|=5$, the boundary curves are 10 exceptional curves on the del Pezzo surface 
$\overline{M}_{0,S}$ corresponding to 10 different distinguished partitions of $S$
whose components have cardinalities $(1,1,1,2)$. They define 10 different
Chow classes.

\medskip

{\bf 3.7. Example: $\overline{M}_{0,6}.$} There are two combinatorial types of unlabeled trees $\tau$
corresponding to boundary curves. Below we draw their subgraphs consisting of all vertices and edges,
and mark them with the numbers of tails at each vertex.
$$
3\bullet -\bullet 1 - \bullet 2 \quad \quad\quad  2\bullet -\bullet 2 -\bullet 2  
$$

If we take into account possible labellings by $S$, we will get  60 boundary curves of the first type and 45 boundary curves
of the second type. They form two different $S_6$--orbits.

\medskip

If $\tau$ is of the first type, then $c(\sigma )=8$ for all 3 partitions $\sigma \in P(\tau )$.
 The set $N(\tau )$ contains unique partition $\sigma $,
with $c(\sigma )=9$. Applying Proposition 3.4, we get
$$
(-K_6, C_{\tau})=1.
$$

\medskip

If $\tau$ is of the second type,  we have respectively $c(\sigma )=8,9,9$ for
$\sigma\in P(\tau )$. The set $N(\tau )$ consists of 2 partitions $\sigma$,
with $c(\sigma )=8$. Applying Proposition 3.4, we get
$$
(-K_6, C_{\tau})=0.
$$

\medskip

Chow classes of the boundary curves for $n=6$ are extremal rays of the Mori cone.
There are 20 classes of the first type and 45 classes of the second type.

\medskip

{\bf 3.8. Example: $\overline{M}_{0,7}.$} Similarly, there are four combinatorial types of unlabeled trees $\tau$
corresponding to boundary curves. 
$$
A:\quad 3\bullet -\bullet 1 - \bullet 1 -\bullet 2  \quad \quad\quad B:\quad  2\bullet -\bullet 2 - \bullet 1 -\bullet 2
$$
and
$$
C:\ 3\bullet -\bullet^{/^{\bullet 2}}_{\setminus_{\bullet 2}} \quad \quad\quad D:\quad  2\bullet -1\bullet^{/^{\bullet 2}}_{\setminus_{\bullet 2}}
$$
\bigskip

Here the numerology looks as follows.

\bigskip

{\it Type A.}  We have $ c(\sigma )=10$ for all $\sigma\in P(\tau )$;
$|N(\tau )|=1, c(\sigma )=12$ for  $\sigma\in N(\tau )$. Hence
$$
(-K_7, C_{\tau})=1.
$$
Finally, there are 420 labeled trees/boundary curves of this type.

\bigskip

{\it Type B.}  We have $c(\sigma )=10, 12, 12$  for  $\sigma\in P(\tau )$;
$|N(\tau )|=2, c(\sigma )=10, 12$ for  $\sigma\in N(\tau )$. Hence
$$
(-K_7, C_{\tau})=0.
$$
There are 630 boundary curves of this type.

\bigskip

{\it Type C.}  We have $ c(\sigma )=10$ for all $\sigma\in P(\tau )$;
$|N(\tau )|=1, c(\sigma )=12$ for  $\sigma\in N(\tau )$. Hence
$$
(-K_7, C_{\tau})=1.
$$
There are 105 boundary curves of this type.

\bigskip

{\it Type D.}  Finally, here $ c(\sigma )=12$ for all $\sigma\in P(\tau )$;
$|N(\tau )|=3, c(\sigma )$=10 for  $\sigma\in N(\tau )$, and
$$
(-K_7, C_{\tau})=-1.
$$
There are 105 boundary curves of this type.

\bigskip

In the Chow group, there are 35 classes of types A and C altogether, 210 classes of type B,
and 105 classes of type D.

\newpage


\centerline{\bf 4. Gromov--Witten correspondences}

\smallskip
\centerline{\bf  for boundary curves in $\overline{M}_{0,S}$.}

\setcounter{section}{4}

\medskip

In this section we  will state and prove the main theorem of this paper.   We start with some preparation.

\medskip

{\bf 4.1. Preparation: combinatorics.} In this section, we choose and fix two disjoint finite sets 
$S$ and $\Sigma$.
Assume that $|S|\ge 4$, $|\Sigma |\ge 3$. 

\smallskip

Fix one element $s_0\in S$. Choose and fix a  distinguished partition $\Pi$ of $S$ into four disjoint nonempty subsets (cf. 3.2 above). Denote by $S(\Pi)$ the
set, elements of which are components  of $\Pi$. Thus, $|S(\Pi)|=4$.
Denote by $\bullet \in S(\Pi )$ the component of $\Pi$ that contains the marked element $s_0\in S$.

\smallskip

The sets $P(\Pi )$ and $N(\Pi )$ are defined as in 3.2.1. In our setup, the three--element set $P(\Pi )$
is canonically bijective to two more sets:

\smallskip

a) The set of stable unordered  partitions of $S(\Pi )$ into two parts (each consisting of two elements).

\smallskip

b) The set $S(\Pi )\setminus \{\bullet\}$: any $j\in S(\Pi )\setminus \{\bullet\}$ corresponds to the
partition $S(\Pi )= (\{\bullet ,j\}\sqcup S(\Pi )\setminus \{\bullet ,j\})$. We have already used this trick in sec. 2.1,
and here we will use it again  transporting the results of sec. 2 to a new context.

\smallskip

Slightly abusing notation, we will sometimes consider these last identifications as identical maps.

\smallskip

Being more fussy, we can say that our constructions are functorial on the category of pointed finite sets $S$ 
with bijections. Eventually, they must be extended to the category of marked trees 
(and more general modular graphs) encoding boundary combinatorial types of curves and maps.
Dependence of our geometric construction on the target boundary curve class $\beta$ is reflected in the dependence
of its combinatorial side on $\Pi$.

\medskip

{\bf 4.2. Preparation: geometry.} We intend to show that  results of sec. 1.5--1.6
are applicable in the present situation.

\smallskip

More precisely, specialize the objects, introduced in 1.5  in the following way:
\begin{align}\label{Eq.: W,E,b,beta}
W:= \overline{M}_{0,S}, \quad E:= \overline{M}_{\Pi}, \quad b:= b_{\Pi},\quad  \beta := \beta (\Pi )
\end{align}
(cf. (3.1)).

\smallskip

Furthermore, specialize the objects described in 1.6:
\begin{align}\label{Eq.: B,C,p}
B:= B_{\Pi},\quad C:= \overline{M}_{0, F_{\pi}(v_0)},\quad p:= p_{\Pi}
\end{align}
(cf. (3.2), (3.3)).

\medskip

{\bf 4.3. Proposition.}  {\it The assumptions 1.6 (a)--(d) hold for \eqref{Eq.: W,E,b,beta}--\eqref{Eq.: B,C,p}.}

\medskip

{\bf 4.4. Proof of Proposition 4.3.}  The  assumptions 1.6 (a) and (b) hold by definition. 

\smallskip

{\bf 4.4.1. Assumption (c).} Let $C$ be a closed fiber of $p:\, \overline{M}_{\Pi}\to B_{\Pi}$. We have already used the fact that it is isomorphic to $ \textbf{P}^1$. Let $j:\, C\to \overline{M}_{0,S}$ be the natural closed  embedding. We assert that
\begin{align}\label{Eq.: Direct Sum Decomposition of the Tangent}
j^*(\mathcal{T}_{\overline{M}_{0,S}}) \cong \mathcal{O}(2) \oplus \mathcal{O}^{n-4-|N(\Pi )|} \oplus \mathcal{O}(-1)^{|N(\Pi )|}
\end{align}
where  $\mathcal{T}_{\overline{M}_{0,S}}$ is the tangent sheaf and $\mathcal{O}:=\mathcal{O}_C$.
 
\smallskip
 
In fact,  consider the embedding $i:\, C\to\overline{M}_{\Pi} $ and  the natural filtration
\begin{align}\label{Eq.: Filtration on the Tangent}
\{0\} \subset \mathcal{T}_C \subset  i^*(\mathcal{T}_{\overline{M}_{\Pi}})\subset  j^*(\mathcal{T}_{\overline{M}_{0,S}}).
\end{align}
The consecutive summands in \eqref{Eq.: Direct Sum Decomposition of the Tangent} correspond to the consecutive quotients of \eqref{Eq.: Filtration on the Tangent}. Namely, $\mathcal{T}_C  \cong \mathcal{O}(2)$;   $i^*(\mathcal{T}_{\overline{M}_{\Pi}})/\mathcal{T}_C$ is trivial of the rank
\begin{align*}
\dim \, B_{\Pi} =\sum_{v\in V_{\pi}} (|F_{\pi}(v)|-2) = |S|-4-|N(\Pi )|,
\end{align*}
finally, the last isomorphism follows from (3.1) and (3.4).
 
\smallskip
 
From \eqref{Eq.: Direct Sum Decomposition of the Tangent} we see that $H^1(C,j^*(\mathcal{T}_{\overline{M}_{0,S}}))=0$.

\smallskip

{\bf 4.4.2. Assumption (d): preparation I.} Any curve $X$ in $\overline{M}_{0,S}$ of class $\beta(\Pi)$ is a closed fiber of $p_{\Pi} \colon \overline{M}_{\Pi}\to B_{\Pi}$. In fact, by (3.1) and (3.4) the curve $X$ is contained in $\overline{M}_{\Pi}$. Below we will show that it is indeed a fiber of $p_{\Pi}$ by analysing degeneration patterns of fibers of the universal family $\overline{C}_{0,S}$ over points of $X$.

\smallskip

Let $\sigma$ be the dual graph of the curve from the universal family  $\overline{C}_{0,S}$ over a generic point of $X$. We know that $\sigma$ admits a contraction onto $\pi$. If $p_{\Pi}(X)$ is not a point, then $X$ must contain a point over which the dual graph $\sigma^{\prime}$ of the universal family is not isomorphic to $\sigma$. In this case it must admit a non--trivial contraction $\sigma^{\prime}\to\sigma$. Compose it with the canonical contraction  $\sigma \to \pi$.

\medskip

One of the following two alternatives must hold:

\medskip

{\it (A) There is an edge of $\sigma^{\prime}$ that contracts onto one of the vertices $v\ne v_0$ of $\pi$}.

\smallskip

{\it (B) No edge of $\sigma^{\prime}$ contracts onto one of the vertices $v\ne v_0$, but there is an edge contracting to $v_0$.}

\medskip

Consider the stable 2--partition $\rho$ of $S$ corresponding to the contracting edge, and the respective boundary divisor $D_{\rho}$ in $\overline{M}_{0,S}$. Geometrically, our  assumption (A) implies that $X$ is a curve that does not lie in $D_{\rho}$ but intersects $D_{\rho}$, hence we must have
$$
(D_{\rho}, \beta ) = (D_{\rho},  [X]) >0.
$$
But from (3.4) it follows that if $\rho$ contracts onto a vertex $v\ne v_0$, then  $(D_{\rho}, \beta )=0$. Hence this possibility is excluded.

\smallskip

Consider now the alternative (B). Then we must have $\rho\in P(\Pi )$. This implies the following degeneration pattern of the induced family of curves parametrized by $X$. At a generic point, the tree of the curve consists of one irreducible component $C$ to which trees are attached at $|N(\tau )|$ different points of this component. When the degeneration at a point of $D_{\rho}$ occurs, $C$ breaks down into two components, say, $C_1$ and $C_2$, and the attached trees are distributed among them: some become attached to $C_1$, and remaining ones to $C_2$. What is important here, is that the {\it labeled combinatorial type of each of the attached  trees does not change} -- otherwise we could have used the option (A) which was already excluded. 

\smallskip

But in this case the image $p_{\Pi}(X)$ must land in the product of the {\it open strata} $\prod_{v\ne v_0} M_{F_{\pi}(v)}$. This is possible only if this image is a point because such a product is an affine scheme.

\smallskip

{\bf 4.4.3. Assumption (d): preparation II.} Let $Hilb(\overline{M}_{0,S})$ be the Hilbert scheme of $\overline{M}_{0,S}$. As usual, it can be written as a disjoint union $Hilb(X)=\coprod_P Hilb^P(X)$, where $P$ is a Hilbert polynomial (for this one should fix an ample line bundle) and each $Hilb^P(\overline{M}_{0,S})$ is a quasi-projective scheme.

\smallskip

Let $C$ be a closed fiber of $p_{\Pi} \colon \overline{M}_{\Pi} \to B_{\Pi}$. It defines a closed point $P$ in $Hilb(\overline{M}_{0,S})$. The tangent space to $Hilb(\overline{M}_{0,S})$ at the point $P$ is identified with $H^0(C, N_{C/\overline{M}_{0,S}})$, and the obstruction space with $H^1(C, N_{C/\overline{M}_{0,S}})$. From computations in Section 4.4.1 it follows that 
\begin{align*}
&\dim H^1(C, N_{C/\overline{M}_{0,S}})=0.
\end{align*}
Therefore, $P$ is a smooth point of $Hilb(\overline{M}_{0,S})$. 

\smallskip

Consider the locus in $Hilb(\overline{M}_{0,S})$ parametrizing closed fibers of $p_{\Pi} \colon \overline{M}_{\Pi} \to B_{\Pi}$. It is a connected component of $Hilb(\overline{M}_{0,S})$. Denote it $Y$, and let $U \to Y$ be the universal family over it. We have just seen that $Y$ is a smooth (in particular reduced and irreducible) scheme. Its dimension 
\begin{align*}
&\dim Y= \dim H^0(C, N_{C/\overline{M}_{0,S}})=\dim B_{\Pi},
\end{align*}
which follows from computations in Section 4.4.1.

Therefore, we can identify the universal family $U \to Y$ with the projection $p_{\Pi} \colon \overline{M}_{\Pi} \to B_{\Pi}$.

\smallskip

{\bf 4.4.4. Assumption (d): proof.} We need to show that the canonical morphism of stacks
\begin{align}\label{Eq.: Induced Iso}
\widetilde{b}_{\Pi } \colon  \overline{M}_{0,\Sigma}(\overline{M}_{\Pi},\beta_{\Pi} )  \to \overline{M}_{0,\Sigma}(\overline{M}_{0,S},
\beta (\Pi )),
\end{align}
induced by $b_{\Pi} \colon \overline{M}_{\Pi} \to \overline{M}_{0,S}$ is an isomorphism. Here $\beta_{\Pi}$ is the Chow class of a  fiber of  $p_{\Pi} \colon \overline{M}_{\Pi}\to B_{\Pi}$.

\smallskip

One $T$--point of  $ \overline{M}_{0,\Sigma}(\overline{M}_{0,S},\beta (\Pi ))$ is a family of connected prestable curves $p_T \colon \mathcal{C}_T \to T$ together with a stable map $f_T$ of the class $\beta (\Pi )$ and labeled sections 
\begin{align*}
f_T:\,\mathcal{C}_T\to \overline{M}_{0,S},\quad  x_{j, T}:\, T\to\mathcal{C}_T, \ j\in \Sigma .
\end{align*}
Below we will show that any such map $f_T$ can be factored through $b_{\Pi}$. Since $b_{\Pi}$ is a closed embedding, such factorization is unique if it exists.

\smallskip

Consider the diagram
\begin{align*}
\xymatrix{
C_T \ar[rr]^{f_T \times p_T} \ar[dr]_{p_T}  & & \overline{M}_{0,S} \times T \ar[dl]^{pr_T}\\
&T
}
\end{align*}
provided by the $T$--point. Since $\beta (\Pi )$ is indecomposable (Proposition 3.5), for any geometric fiber of $\mathcal{C}_T/T$, $f_T$ must contract to a point each its component excepting one. On the uncontracted component it is a closed embedding.

\smallskip

{\bf Irreducible geometric fibers.} Assume that all geometric fibers of $p_T$ are irreducible and hence $f_T\times p_T$ induces closed embeddings on all geometric fibers. By faithfully flat descent it is then a closed embedding on all fibers. Therefore, the fiber of $f_T\times p_T$ at a point $s \in \overline{M}_{0,S} \times T$ is either empty or $\kappa(s)$-isomorphic to $\text{Spec}(\kappa(s))$, where $\kappa(s)$ is the residue field at $s$.

Since $p_T$ and $pr_T$ are proper, the morphism $f_T\times p_T$ is also proper. By [EGA], Proposition 8.11.5 it implies that $f_T\times p_T$ is a closed embedding. Thus, we see that if we forget the sections $(x_{j,T})$ the stable morphism $(C_T, (x_{j,T}), f_{T})$ gives us a $T$--point of the Hilbert scheme of $\overline{M}_{0,S}$. 

Therefore, by Section 4.4.3 the diagram
\begin{align*}
\xymatrix{
C_T \ar[r]^{f_T} \ar[d]  & \overline{M}_{0,S}\\
T
}
\end{align*}
is obtained from 
\begin{align*}
\xymatrix{
\overline{M}_{\Pi} \ar[r]^{b_{\Pi}} \ar[d]_{p_{\Pi}}  & \overline{M}_{0,S}\\
B_{\Pi}
}
\end{align*}
by a unique pullback. Hence, the stable map $(C_T, (x_{j,T}), f_{T})$ factors through $\overline{M}_{\Pi}$.

\medskip

{\bf General case.} Let $(C_T, (x_{j,T})_{j \in \Sigma}, f_{T})$ be an arbitrary $\Sigma$-labelled stable map to $\overline{M}_{0,S}$ of class $\beta(\Pi)$, and let $\Sigma' \subset \Sigma$ be the subset that labels sections that land on the non-contracted component of geometric fibers. Consider the induced \textit{prestable} map $(C_T, (x_{j,T})_{j \in \Sigma'	}, f_{T})$. Stabilizing it we get a stable map $(\widetilde{C}_T, (y_{j,T})_{j \in \Sigma'	}, g_{T})$ to $\overline{M}_{0,S}$ of class $\beta(\Pi)$, such that $f_T=g_T \circ st$. In other words, we get a diagram
\begin{align*}
\xymatrix{
C_T \ar[d] \ar[r]^{st} & \widetilde{C}_T  \ar[r]^{g_T}\ar[ld]   & \overline{M}_{0, S}  \\
T
}
\end{align*}
where $\widetilde{C}_T \to T$ has irreducible geometric fibers. As we have seen above, $g_T$ factors through the embedding $b_{\Pi} \colon \overline{M}_{\Pi} \to \overline{M}_{0, S}$, and hence so does $f_T$. 

\medskip

We have shown that any family of stable maps to $\overline{M}_{0,S}$ of class $\beta(\Pi)$ can be factorized uniquely via the closed embedding $b_{\Pi} \colon \overline{M}_{\Pi}  \to \overline{M}_{0,S}$. Therefore, it gives a family of stable maps to $\overline{M}_{\Pi}$ of class $\beta_{\Pi}$. 

This procedure gives a map of $T$--points of stacks appearing in \eqref{Eq.: Induced Iso} for any $T$. One can check that it naturally extends to morphisms of $T$--points and gives a functor 
\begin{align}\label{Eq.: Inverse Functor}
\overline{M}_{0,\Sigma}(\overline{M}_{0,S},\beta (\Pi )) \to \overline{M}_{0,\Sigma}(\overline{M}_{\Pi},\beta_{\Pi} ). 
\end{align}
Moreover, one can easily check that \eqref{Eq.: Inverse Functor} is indeed inverse to \eqref{Eq.: Induced Iso}. We leave these checks to the reader.

\medskip

{\bf 4.5. The final summary.} We will now briefly  restate the results of stepwise calculations 
of sec. 1 and 2 in our current situation \eqref{Eq.: W,E,b,beta} -- \eqref{Eq.: B,C,p}.
\medskip

{\bf 4.5.1. Step I: Gromov--Witten correspondences for the target space $\overline{M}_{0,S(\Pi)}$.}
We reproduce here the main result of sec. 2 applied to the target space $\overline{M}_{0,S(\Pi)}$ 
and its fundamental class $\textbf{1}$. Notice that the sets denoted $S$ (resp. $P$) in sec. 2
are now $S(\Pi )$ (resp. $P(\Pi )$), and $S(\Pi )= P(\Pi)\sqcup \{\bullet \}$.

\smallskip

According to   Proposition 2.3, we have:
\begin{align}\notag
I_{0,\Sigma}(\overline{M}_{0,P(\Pi )\sqcup \{\bullet \}},\textbf{1} )= (\overline{ev}, \psi_{P(\Pi )})_*
 ([\overline{M}_{0,\Sigma\sqcup P(\Pi )}])\in \\ \label{Eq.: GW Corr. M_0,4 beta=1}
\in  A_{|\Sigma |}( (\overline{M}_{0,P(\Pi )\sqcup \{\bullet \}})^{\Sigma}
\times \overline{M}_{0,\Sigma}    ).
\end{align}

\medskip

{\bf 4.5.2. Step II: Gromov--Witten correspondences for the target space $B_{\Pi}$ and zero beta--class.} According to the Example 1.3, we have:
\begin{align}\label{Eq.: GW Corr. B_Pi beta=0}
I_{0,\Sigma}(B_{\Pi}, 0 ) = [\Delta_{\Sigma}(B_{\Pi})\times \overline{M}_{0,\Sigma}]\in A_*(B_{\Pi}^{\Sigma}\times \overline{M}_{0,\Sigma}).
\end{align}
Here $\Delta_{\Sigma}(B_{\Pi})$ is the diagonal in the cartesian product  $B_{\Pi}^{\Sigma}$ of   
$\Sigma$ copies of $B_{\Pi}$.
\medskip

{\bf 4.5.3. Step III: Gromov--Witten correspondences for the target space $\overline{M}_{\Pi}$ and fiber  beta--class.}  In this subsection,  $\beta_{\Pi}$ is the Chow class of a fiber of the projection $\overline{M}_{\Pi}\to B_{\Pi}$. We now have a canonical splitting
\begin{align*}
\overline{M}_{\Pi}= B_{\Pi}\times  \overline{M}_{0,P(\Pi )\sqcup \{\bullet \}}
\end{align*}
since $F_{\pi}(v_0)$ is identified with $S(\Pi )=P(\pi )\sqcup \{\bullet\}$ (cf. 3.3).

\medskip

Thus using \eqref{Eq.: GW Corr. M_0,4 beta=1} and \eqref{Eq.: GW Corr. B_Pi beta=0},  we have
\begin{align}\label{Eq.: GW Corr for M_Pi}
I_{0,\Sigma}(\overline{M}_{\Pi},\beta_{\Pi} )= \widetilde{\Delta}^!([\Delta_{\Sigma}(B_{\Pi})\times \overline{M}_{0,\Sigma}]\otimes (\overline{ev}, \psi_{P(\Pi )})_* ([\overline{M}_{0,\Sigma\sqcup P(\Pi )}])) \, .
\end{align}

\bigskip

To summarize, we have proved our final theorem, a specialization of Proposition 1.6.1:

\medskip

{\bf 4.6. Theorem.} {\it  The structure embedding $b_{\Pi}:\, \overline{M}_{\Pi} \to \overline{M}_{0,S}$ induces a canonical isomorphism
\begin{align*}
\widetilde{b}_{\Pi}:\  \overline{M}_{0,\Sigma}(\overline{M}_{\Pi},\beta_{\Pi} )  \to \overline{M}_{0,\Sigma}( \overline{M}_{0,S}, \beta (\Pi )).
\end{align*}
where $\beta_{\Pi}$ is the Chow class of a  fiber of  $p_{\Pi}:\, \overline{M}_{\Pi}\to B_{\Pi}.$
\smallskip
This isomorphism $\widetilde{b}_{\Pi}$ is compatible with evaluation/stabilization morphisms for both
moduli spaces and induces the identity
\begin{align}\label{Eq.: Corr. Direct Image}
I_{0,\Sigma}(\overline{M}_{0,S},\beta (\Pi )) = (b_{\Pi}^{\Sigma}\times id)_*( I_{0,\Sigma}(\overline{M}_{\Pi},\beta_{\Pi}))
\end{align}
where
$$
b_{\Pi}^{\Sigma}\times id :\quad  \overline{M}_{\Pi}^{\Sigma}\times  \overline{M}_{0,\Sigma} \to (\overline{M}_{0,S})^{\Sigma}\times  \overline{M}_{0,\Sigma} .
$$

\smallskip

The rhs of \eqref{Eq.: Corr. Direct Image} is given by \eqref{Eq.: GW Corr for M_Pi}.}

\bigskip

{\bf 4.7. Gromov--Witten numbers.} In this subsection, we will specialize formula (1.18)
to our situation in order to calculate numerical invariants of Chow classes of boundary curves.

\smallskip
Let $\gamma_j\in H^{2d_j}(\overline{M}_{0,S})$ be a family of cohomology classes
indexed by $j\in \Sigma$. If  $\sum_{j\in\Sigma}d_j=\dim \,B_{\Pi}$, then
the correspondence 
$$
I_{0,\Sigma}(\overline{M}_{0,S},\beta (\Pi))\in A_*((\overline{M}_{0,S})^{\Sigma}\times
\overline{M}_{0, \Sigma})
$$ 
maps
$\otimes_{j\in\Sigma}\gamma_j \in (H^*(\overline{M}_{0,S}))^{\otimes \Sigma}$
to a class of maximal dimension in $H^*(\overline{M}_{0,\Sigma}).$ The degree
of this class is denoted
$$
\langle  I_{0,\Sigma,\beta(\Pi )}^{\overline{M}_{0,S}}\rangle (\otimes_{j\in\Sigma}\gamma_j).
$$
Generally, this degree is  the virtual number of  stable maps of pointed curves of class $\beta (\Pi )$
satisfying the incidence conditions $f(x_j)\in \Gamma_j$, where $(\Gamma_j)$ are cycles
in general position whose dual classes are $\gamma_j$:
$$
f:\,(C;(x_j\,|\,j\in\Sigma)) 
\to \overline{M}_{0,S} ,
$$
whenever such incidence conditions are strong enough to enforce existence only of finite
(virtual) number
of such maps.  In our unobstructed case, this virtual number is the actual 
number of such maps whenever the incidence cycles are in general position.
\smallskip
Recall also that this number is polylinear in $(\gamma_j)$.
\medskip

{\bf 4.7.1. Proposition.} {\it We have
\begin{align}\label{Eq.: GW numbers}
\langle  I_{0,\Sigma,\beta(\Pi )}^{\overline{M}_{0,S}}\rangle (\otimes_{j\in\Sigma}\gamma_j)
=\text{deg}\left( \cap_{j\in\Sigma} pr_{B_{\Pi}*}\circ b_{\Pi}^*(\gamma_j) \right) \, .
\end{align}
}
\smallskip
{\bf Sketch of proof.} Skipping a clumsy but straightforward formal derivation of \eqref{Eq.: GW numbers} from (1.18), we describe the geometric content of this counting formula in the general situation axiomatized in 1.6.

\smallskip

First of all, (1.18) reduces the count to the case of an incidence condition represented by some cycles in $E=\overline{M}_{\Pi}$: in fact,  $b_{\Pi}^*(\gamma_j)$ are  represented by $\Gamma_j\cap \overline{M}_{\Pi}$ in the case of transversal intersections.

\smallskip

Now, in $\overline{M}_{\Pi}$ the incidence cycles can be replaced by ones of the form $\Delta_j\times c_j +\Delta^{\prime}_j\times C$ where $c_j$ are points on a projective line $C$ as in \eqref{Eq.: B,C,p} corresponding to the decomposition $\overline{M}_{\Pi}=B_{\Pi}\times C$.

\smallskip

Assume first that $\Delta_j^{\prime}\ne 0$ for some $j=j_0$. If for such an incidence condition there is a fiber $C_0$ of $\overline{M}_{\Pi} \to B_{\Pi}$ satisfying it at all, then the number of relevant pointed stable maps must be infinite, because   $x_{j_0}$ can be chosen arbitrarily along this fiber. Hence decomposable cycles containing at least one factor of the form $\Delta^{\prime}_j\times C$ give zero contributions to \eqref{Eq.: GW numbers}.

\smallskip

Now consider the case of incidence conditions of the form $\Delta_j\times c_j$ for all $j\in\Sigma$. Let $\Delta_j=pr_{B_{\Pi}}(\Delta_j\times c_j)$ be in a general position in $B_{\Pi}$ so that the intersection cycle $\cap_{j\in\Sigma} \Delta_j$ is a sum of points $y_a\in B_{\Pi}$, of multiplicity one each. We can also lift $\Delta_j$ arbitrarily to $\overline{M}_{\Pi}$, that is choose all $c_j\in C$ pairwise distinct, and consider $\Delta_j\times c_j$ as a geometric incidence condition representing the initial cohomological incidence condition $(\gamma_j)$. 

\smallskip

After that the geometric count becomes straightforward: each point $y_a$ produces one fiber of the class $\beta (\Pi)$ intersecting each $\Delta_j\times c_j$ at one point corresponding to $c_j$. 

\smallskip

The number of $(y_a)$ is the right hand side of \eqref{Eq.: GW numbers}, and the curve count interprets the left hand side of \eqref{Eq.: GW numbers}.

\bigskip

\centerline{\bf 5. Examples and remarks}

\medskip

{\bf 5.1. ``Naturality'' of Gromov--Witten correspondences.} 
In this subsection we try to make somewhat more precise our guess 0.2.1.
To this end we recall first, that natural objects in the relevant category are
moduli spaces $\overline{M}_{\tau}$, and natural morphisms/correspondences
are those ones that  are produced from morphisms in the category
of modular graphs. The latter include contractions, forgetful morphisms,
relabeling morphisms etc., cf. [BehMa].

\smallskip

The least controllable characteristic of GW--correspondences is their dependence
on the argument $\beta$ in the relevant Mori cone. So far we
have considered only boundary $\beta$'s, and they are, of course,
``natural'' by definition.
\smallskip

In this subsection we will show  that,
keeping notation of section 4, we may naturally encode
most of the relevant combinatorial and geometric information in one
moduli space 
$\overline{M}_{0, \Sigma \times (S\setminus \{s_0\})}$
and a configuration of certain of its boundary strata. This is only a tentative suggestion, we do not develop it
fully, because we still lack even a conjectural description of the situation for more
general $\beta$'s.

\medskip

{\bf 5.1.1. The tree $\textbf{T}$.}  The tree $\textbf{T}$ has one special vertex that will
be called {\it central} one and denoted $v_c$. Its flags are bijectively labeled by $\Sigma$:
we will use the notation
$$
F_{\textbf{T}}(v_c):= \{\langle j\rangle \,|\, j\in \Sigma \} \ .
\eqno(5.1)
$$
The remaining vertices constitute a set bijective to  $\Sigma\times \{s_0\}$. Together with (5.1),
this bijection is a part of structure, and we may refer to a non--central vertex $v\in V_{\textbf{T}}$ as
$v_j:=\langle j,s_0\rangle ,\ j\in \Sigma$.
\smallskip

Furthermore, we put
$$
F_{\textbf{T}}(v_j ):= \{j\} \times S  = \{(j,s)\,|\,s\in S\}\ .
\eqno(5.2)
$$

\smallskip

Thus, the standard  identification of $\overline{M}_{\textbf{T}}$ with the product of moduli spaces 
$ \prod_{v\in V_{\textbf{T}}}\overline{M}_{0,F_{\textbf{T}}(v)}$ corresponding to stars of all vertices,
can be rewritten as
$$
\overline{M}_{\textbf{T}} = (\overline{M}_{0,S})^{\Sigma}\times \overline{M}_{0,\Sigma}
\eqno(5.3)
$$
where the last factor corresponds to the central vertex.

\smallskip

{\it (ii) Edges.} The  flag $\langle j \rangle$ attached to the central vertex  (see (5.1)) is grafted to the flag 
$(j,s_0)$
incident to the vertex $v_j$ (see (5.2)).  There are no more edges.

\smallskip

Thus, the central vertex
carries no tails, and the set of edges $E_{\textbf{T}}$ is naturally bijective to $\Sigma$.
The set of tails is
$$
T_{\textbf{T}}=\coprod_{j\in \Sigma} (F_{\textbf{T}}(\langle j,s_0\rangle )\setminus (j,s_0))=
\coprod_{j\in \Sigma} (\{j\}\times (S\setminus \{s_0\}))   \cong
 \Sigma \times (S\setminus \{s_0\}) .
 \eqno(5.4)
$$
\smallskip
If we interpret the last set in (5.4) as the set of labels of tails , then
the described above set of edges of $\textbf{T}$ determines the canonical embedding
of $\overline{M}_{\textbf{T}}$ as a boundary stratum:
$$
\overline{M}_{\textbf{T}} \hookrightarrow   \overline{M}_{0, \Sigma \times (S\setminus \{s_0\}) }
\eqno(5.5)
$$
This embedding corresponds to full contraction of all edges  of $\textbf{T}$ to the star with
flags $T_{\tau}$.

\smallskip

We will now encode information about $\Pi$ into another tree $\textbf{T}(\Pi )$,
together with its contraction onto $\textbf{T}$.
 
 \medskip
 
{\bf 5.1.2.  The tree $\textbf{T}(\Pi )$.}  Briefly, to get $\textbf{T}(\Pi )$, we replace
each non--central vertex $v_j,$ $j\in \Sigma$, by a copy $\pi_j$ of the tree $\pi$ described in 3.2.

\smallskip

More precisely, the special vertex of $\pi_j$ denoted $v_{0,j}$ now carries tails (5.2)
distributed among other vertices of $\pi_j$ according to $\Pi$, and its tail $(j,s_0)$ is grafted in $\textbf{T}(\Pi )$
to the same flag $\langle j\rangle$ of its central vertex as it was in $\textbf{T}$.

\smallskip

The contraction $\textbf{T}(\Pi )\to \textbf{T}$ contracts each $\pi_j$ to the star of $v_j$,
and is identical on the stars of the central vertices. Combining the relevant boundary
morphism with (5.5), we get the diagram of strata embedding

$$
\overline{M}_{\textbf{T}(\Pi)} \hookrightarrow \overline{M}_{\textbf{T}} \hookrightarrow   \overline{M}_{0, \Sigma \times (S\setminus \{s_0\}) }
\eqno(5.5)
$$
The intermediate and final correspondences considered in sec. 4, can be expressed
using the geometry of (5.5).

\medskip

{\bf 5.2. Using the Reconstruction Theorems.}  For a general target $W$,
if the Chow ring $A^*(W)$ (with coefficients in $\textbf{Q}$) coincides with the whole $H^*(W)$ and
 is generated by $A^1(W)$, then the total motivic quantum cohomology of $W$ of genus zero understood as the
 family of $I$--correspondences  is completely determined
by triple correlators (3--point GW--invariants) of codimension zero. This follows from the First and
the Second Reconstruction Theorems of [KoMa1].

\smallskip

In any case, these triple correlators are precisely coefficients of small quantum cohomology
as a formal series in $q^{\beta}$. Hence in the same assumptions the total quantum cohomology
is completely determined by the small quantum multiplication in $H^*(V)[[q^{\beta}]]$:
$$
\Delta_a\cdot \Delta_b= \Delta_a \cup \Delta_b + \sum_{\beta\ne 0} \sum_{c\ne 0}
\langle\Delta_a\Delta_b\Delta_c\rangle_{\beta}\Delta^c q^{\beta} .
$$
Here $(\Delta_a)$ is a basis of $H^*$ such that $\Delta_0$ is identity, 
$g_{ab}= (\Delta_a,\Delta_b)$, $(g^{ab})$ the inverse matrix to $(g_{ab})$, and 
and
$\Delta^a:=\sum_b g^{ab}\Delta_b$. 

\smallskip

This is applicable to all $\overline{M}_{0,n}$.

\smallskip

In turn, the associativity equations allow one to express all triple  correlators through
a part of them. We will now make explicit this subset   for $\overline{M}_{0,n}$. 

\medskip

{\bf 5.3. A generating subset of triple correlators.}  Put $|\Delta | =i$ for $\Delta \in H^{2i}( \overline{M}_{0,n}).$
(No confusion with cardinality $|S|$ of a set $S$ should arise).
Then all invariants can be recursively calculated through 3--point invariants $\langle \Delta_a\Delta_b\Delta_c\rangle_{\beta}$
with $\Delta_c$ divisorial, $|\Delta_a|, |\Delta_b|\ge 1$, $\beta \ne 0$, and
$$
|\Delta_a|+|\Delta_b|= (-K_n,\beta )+n-4.
$$
where $K_n$ is the canonical class of  $\overline{M}_{0,n}$. Hence, $\beta$ are restricted by
$$
2-(n-3)\le (-K_n,\beta )-1\le n-3.
$$
\smallskip
See [KoMa1], Theorem 3.1, with the following easy complements. If $|\Delta_a|$ or $|\Delta_b|=0$, $\beta \ne 0$,
then the respective GW--invariant is 0 because of [KoMa1], (2.7). If $\beta =0$,
we can use [KoMa1], (2.8). It remains to consider the following list of parameters:
$$
6-n\le (-K_n,\beta ) \le n-2,
$$
$$
2 \le |\Delta_a|+|\Delta_b|= (-K_n,\beta )+n-4\le 2n-6
$$
Finally, if $\Delta$ is a divisorial class with $(\Delta ,\beta )=0$, then
$\langle \Delta '\Delta '' \Delta\rangle_{\beta} =0$ for any 
$\Delta , '\Delta '' $ due to the Divisor Axiom.
\bigskip

{\bf 5.3.1. Tables for the first values of $n$}.

\vskip0,5cm

$(-K_5,\beta)\quad |$\ \quad 1\ \quad  \  2\ \quad  \quad  3 

{\bf ------------------------------------------}

$(|\Delta_a|, |\Delta_b|)|$\  (1,1)\ \    (1,2)  \  (2,2)

\vskip0,5cm

$(-K_6,\beta)\quad |$\ \quad 0 \ \quad 1\ \quad\quad  2\ \quad\quad  \  \ 3\ \quad  \quad \quad   4 

{\bf ------------------------------------------------------}

$(|\Delta_a|, |\Delta_b|)|$\  (1,1)\ \    (1,2) \  \  (2,2)\ \quad (2,3) \quad  (3,3)

\quad\quad\quad\quad\quad\quad\quad\quad\quad\quad\quad \ (1,3)

\vskip0.5cm

$(-K_7,\beta)\quad |$\  \quad -1   \ \quad 0 \ \quad\quad 1\ \quad\quad \  2\ \quad\quad  \  \ 3\ \quad  \quad \quad   4  \quad \quad \ 5 

{\bf -----------------------------------------------------------------------}

$(|\Delta_a|, |\Delta_b|)|$\quad  (1,1)\ \    (1,2) \  \  (2,2)\ \quad (2,3) \quad  (3,3) \quad (3,4) \quad (4,4)

\quad\quad\quad\quad\quad\quad\quad\quad\quad\quad\quad\quad (1,3) \quad (1,4)\quad \ (2,4)

\bigskip

\bigskip

Notice that  $\overline{M}_{0,5}$ is the del Pezzo surface of degree 5, in particular, its anticanonical class is
ample and hence the generating subset of triple correlators is finite.
In fact, generating sets for  del Pezzo surfaces are collected in [BaMa1]. It is known also that
all del Pezzo surfaces have generically semisimple quantum cohomology,
and more generally, this remains true for blow ups of any finite set of points
on or over $\textbf{P}^2$ (A.~Bayer [Ba]).

\smallskip

Already for $\overline{M}_{0,6}$ the situation is  more mysterious. 
For 45 out of 105 generators of the cone of $\beta$'s we have $(-K_6,\beta )=0$. 
Hence our generating list above is in principle infinite. Semisimplicity is an open question as well.
For $n\ge 7$ the difficulties grow.

\bigskip

{\bf 5.4. Strategies of computation.}  A  possible way to compute some
Gromov--Witten invariants of $\overline{M}_{0,n}$ with non--boundary
$\beta$'s consists in choosing
a birational morphism   $p_n :\,\overline{M}_{0,n}\to X_n$ such that

\medskip

{\it a) (Sufficiently many) GW--invariants of $X_n$ are known/computable.

\smallskip

b) Morphism $p_n$ is such that there exist ``naturality'' formulas that allow one 
to compute (some) GW--invariants of
$\overline{M}_{0,n}$ through (some) GW--invariants of  $X_n$.} 

\medskip

For "naturality" results see, e.~g., [LeLWa], [MauPa], [Hu1], [Hu2], [BrK] (this paper contains
 corrections to [Hu1]),
[Mano1], [Mano2], etc.  We will discuss the relevant classes of morphisms below.
 
\medskip

{\bf 5.4.1.  Blowing $\overline{M}_{0,n}$ down.}
The following choices of morphisms 
seem  promising for application of this strategy, at least for small values of $n$.

\smallskip

(i) $X_n=\textbf{P}^{n-3}$, $p_n$ = Kapranov's morphism, representing
$\overline{M}_{0,n}$ as the result of consecutive blowing up $n-1$ points,
preimages of lines connecting pairs of these points, preimages of planes, passing through
triples of them etc., cf. [HaT].  It involves forgetting
the $n$--th point, then fixing $p_1,\dots ,p_{n-1}\in \textbf{P}^{n-3}$. 

\medskip

(ii) $X_n=(\textbf{P}^1)^{n-3}$, $p_n$ is a similar morphism that was described explicitly
by Tavakol.

\medskip

(iii) $X_n= \overline{L}_{n-2}$, the Losev--Manin moduli space parametrizing stable
chains of $\textbf{P}^1$'s  with marked points and a specific stability condition;
$p_n$ the respective stabilization morphism.

\medskip

It makes sense not just to use $ \overline{L}_{n-2}$ in order to help calculate GW--invariants of   $\overline{M}_{0,n}$,
but to treat these moduli spaces as   replacements of $\overline{M}_{0,n}$ in their own right.
In fact, one can define GW--invariants based upon  $ \overline{L}_{n-2}$,
essentially, no information is lost thereby: see [BaMa2]. 
\medskip
The spaces  $\overline{L}_{n-2}$
are toric, and have the largest Chow ring of these three examples.
These manifolds are
not Fano for $n\ge 6$, but
according to [Ir], any toric manifold has generically semisimple quantum cohomology,
therefore it can be more accessible.

\medskip

(iv) Finally, combining two or more forgetful  morphisms, one can birationally map
$\overline{M}_{0,n}$ and $\overline{L}_{n-2}$ onto products of similar manifolds,
thus opening a way to an inductive calculation of GW--invariants. Here is the simplest 
 example: for $n\ge 5$,  forgetting at first $x_n$, and then all points except for $(x_1,x_2,x_3,x_n)$,
 we get a birational morphism
 $$
\overline{M}_{0,n} \to\overline{M}_{0,n-1} \times  \overline{M}_{0,4},\quad   \overline{M}_{0,4} \cong \textbf{P}^1.
$$
GW--invariants of a product can be calculated via the general quantum K\"unneth formula
whenever they are known for lesser values of $n$.

\smallskip

For our main preoccupation here, that of understanding motivic properties of 
quantum cohomology correspondences, versions of this last
suggestion are most promising.

\bigskip
\centerline{\bf References}

\medskip

[Ba] A.~Bayer. {\it  Semisimple quantum cohomology and blow-ups.}  Int. Math. Res. Not.  2004,  no. 40, 2069"1¤783. 

\smallskip

[BaMa1] A.~Bayer, Yu.~Manin.   {\it  (Semi)simple exercises in quantum cohomology
} In: The Fano Conference Proceedings,
ed. by A. Collino, A. Conte, M. Marchisio, Universit\'a
di Torino, 2004, 143--173. Preprint math.AG/0103164

\smallskip

[BaMa2] A.~Bayer, Yu.~Manin. {\it Stability Conditions, wall-crossing and weighted Gromov--Witten invariants.} Moscow Math Journal, vol. 9, Nr 1 (2009)(Deligne's Festschrift), 3--32. e--print math.AG/0607580

\smallskip

[Beh1]  K.~Behrend. {\it Gromov--Witten invariants in algebraic geometry.}
Inv. Math., 127 (1997), 601--617.

\smallskip

[Beh2]  K.~Behrend. {\it  The product formula for Gromov--Witten invariants.}
J. Algebraic Geom.  8  (1999),  no. 3, 529--541. 

\smallskip

[Beh3] K.~Behrend. {\it Algebraic Gromov--Witten invariants.} In: New trends in algebraic geometry (Warwick, 1996). London Math. Soc. Lecture Note Ser., 264, 
Cambridge Univ. Press, Cambridge, 1999, 19--70.

\smallskip

[BehMa]  K.~Behrend, Yu.~Manin. {\it Stacks of stable maps and Gromov--Witten invariants.}
Duke MJ, 85:1 (1996), 1--60.

\smallskip

[BrK] J.~Bryan, D.~Karp. {\it The closed topological string via
the Cremona trasform.} J. Alg. Geom., 14 (2005), 529--542.

\smallskip

[C] Ana--Maria Castravet.  {\it The Cox ring of  $\overline{M}_{0,6}$.}
Trans. AMS, 361 (2009), Nr. 7, 3851--3878.

\smallskip

[CT]  Ana-Maria Castravet, Jenia Tevelev. {\it Hypertrees, projections, and moduli of
stable rational curves.} arXiv:1004.2553

\smallskip

[EGA] A.~Grothendieck. \textit{\'El\'ements de g\'eom\'etrie alg\'ebrique. IV. \'Etude locale des sch\'emas et des morphismes de sch\'emas. III.} Inst. Hautes \'Etudes Sci. Publ. Math. No. 28 1966 255 pp. 

\smallskip

[FG] G.~Farkas, A.~Gibney. {\it The Mori cones of moduli spaces of pointed curves of small genus.}
Trans. AMS, 355:3 (2003), 1183--1199.

\smallskip

[FuMPh] W.~Fulton, R.~MacPherson. {\it A compactification of configuration spaces.}
Ann. of Math., 130 (1994), 183--225.

\smallskip

[FuPa]  W.~Fulton, R.~Pandharipande. {\it Notes on stable
maps and quantum cohomology.}  In: Proc. Symp. Pure Math., vol. 62, Part 2, 45--96.
Preprint arXiv:9608011

\smallskip

[GM] A.~Gibney, D.~MacLagan. {\it Lower and upper bounds for NEF cones.} arxiv:\
1009.0220

\smallskip

[Hac] P.~Hacking. {\it The moduli space of curves is rigid.} Algebra Number Theory  2  (2008),  no. 7, 809"1¤78. arxiv:0509567
\smallskip

[HaT] B.~Hassett, Yu.~Tschinkel. {\it On the effective cone of the moduli space of pointed
rational curves.} arXIv: 0110231

\smallskip

[Ir] H.~Iritani.  {\it Convergence of quantum cohomology by quantum Lefschetz.}  arXiv:0506236

\smallskip

[Hu1] Jianxun Hu. {\it Gromov--Witten invariants of blow--ups along points and
curves.} Math. Z. 233 (2000), 709--739.

\smallskip

[Hu2] Jianxun Hu. {\it Quantum cohomology of blow--ups of surfaces and
its functoriality property.} Acta Mathematica Scientia,  26B(4),  (2006), 735--743.

\smallskip
[Ka] R.~Kaufmann. {\it The intersection form in $H^*(\overline{M}_{0n})$ and the
explicit K\"unneth formula in quantum cohomology.}
IMRN 19 (1996), 929--952.
\smallskip

[Ke] S.~Keel. {\it Intersection theory of moduli spaces of stable n--pointed curves of genus zero.}
Tr. AMS, 330 (1992), 545--574.

\smallskip

[KeMcK] S.~Keel, J.~McKernan. {\it Contractible extremal rays on  $\overline{M}_{0n}$.}
arxiv: 9607009

\smallskip

[KoMa1]  M.~Kontsevich, Yu.~Manin. {\it Gromov--Witten classes, quantum cohomology, and enumerative geometry.}
Comm. Math. Phys., 164 (1994), 525--562.

\smallskip
[KoMa2]   M.~Kontsevich, Yu.~Manin. {\it Relations between the correlators of the topological
sigma--model coupled to gravity.}
 Comm. Math. Phys., 196 (1998), 385--398.  e--print
alg--geom/970824.

\smallskip
[KoMaKa] M.~Kontsevich, Yu.~Manin (with Appendix by R.~Kaufmann). {\it Quantum cohomology
of a product.} Inv. Math., 124 (1996), 313--339.

\smallskip

[LeLWa] Y.--P.~Lee, H.--W.~Lin, C.--L.~Wang. {\it Flops, motives and invariance
of quantum rings.} Annals of Math.,  172 (2010), 243--290

\smallskip
[Ma] Yu.~Manin. {\it Frobenius manifolds, quantum cohomology, and moduli
spaces.} AMS Colloquium Publications, vol. 47, Providence, RI, 1999,
xiii+303 pp.

\smallskip

[Mano1]  Ch.~Manolache. {\it Virtual pull--backs.} arXiv:0805.2065

\smallskip

[Mano2]  Ch.~Manolache. {\it Virtual push--forwards.} arXiv:1010.2704

\smallskip

[MauPa] D.~Maulik, R.~Pandharipande. {\it A topological view of Gromov--Witten
theory.} Topology, 45 (2006), 887--918. arXiv: 0412503

\end{document}